\documentclass[11pt]{article}
\usepackage[T1]{fontenc}
\usepackage{amsmath,amssymb,amsthm,graphicx}
\usepackage[margin=1in]{geometry}
\usepackage[hidelinks]{hyperref}
\numberwithin{equation}{section}
\newtheorem{proposition}{Proposition}[section]
\newtheorem{theorem}[proposition]{Theorem}
\newtheorem{lemma}[proposition]{Lemma}
\newtheorem{corollary}[proposition]{Corollary}

\theoremstyle{definition}
\newtheorem{example}{Example}[section]
\newtheorem{remark}[proposition]{Remark}
\newtheorem*{remark*}{Remark}
\newtheorem*{definition*}{Definition}
\title{Plane Curves Whose Curvature Is a Function of the Distance to the Tangent Line: Relations Between Karpenkov's Areal and Angular Densities}
\author{Esther Notik\\[2pt]\small Independent researcher, Petah Tikva, Israel\\\small\texttt{esthnt@gmail.com}}
\date{}
\begin{document}
\maketitle
\begin{abstract}
For a plane curve and a point $O$, O. Karpenkov introduced the areal density $A$ and the angular density $B$ as smooth counterparts of the alternating entries of an LLS sequence. Here $A$ is the signed distance from $O$ to the tangent line, and, with the order of the chords in the signed definition fixed here, $A^2B=\kappa$, where $\kappa$ is the signed curvature. A relation between the two densities therefore prescribes the curvature as a function of the distance to the tangent line. Karpenkov asked which curves have $AB$, or $B$ alone, constant. The equation $AB=C$ gives the shrinking, stationary and expanding profiles of planar curve shortening. For constant nonzero $B$, choose the orientation so that $B>0$. Curves with positive support have periodic support or are centred circles. The noncircular closed curves are classified by their turning fractions in $(1/2,1/\sqrt3)$; the centred circle is the only embedded closed solution. Negative support gives smooth open graphs whose two ends approach rays from the origin. On positive support the polar reciprocal is a constant-force orbit, and the original curve is a rotated and scaled velocity curve of that orbit. For $A^3B=C$ we determine which ends extend to infinity, reach a cusp, or spiral into the origin. We describe regularity and the range of turning fractions for the positive-support family $A^pB=C$. Finally, for a general relation $B=\Phi(A)$ the polar reciprocal is an orbit of the central force $L^2/\Phi(1/r)$; we describe the ends of all arcs with $AB<0$ and show that a closed curve with turning fraction $m/n$ has exactly $n(m-1)$ double points. We also describe the curves whose areal or angular density is a linear function of the arclength.
\end{abstract}


\section{Introduction}

The construction of ``infinitesimal continued fractions'' by O.~Karpenkov
\cite{karpenkov,book} associates with an arclength-parametrized $C^2$ plane curve
$\gamma$ and a distinguished point $O$ two functions, the areal density
$A$ and the angular density $B$. Concretely, $A$ is the signed distance
from $O$ to the tangent line, and $B=\kappa/A^2$, where $\kappa$ is the
signed curvature (\S2).

Karpenkov related the areal density $A$ to Kepler's second law
\cite[Proposition~3.2]{karpenkov}. On an arc where $A\ne0$, let a point
move along the curve with speed $v$. The line segment from $O$ to the
moving point sweeps out area at rate $\tfrac12 |A|v$. Thus choosing
$v$ proportional to $1/|A|$ gives equal swept areas in equal times.
The identity $A^2B=\kappa$ (with the sign convention fixed in \S2) relates
the two densities to the Euclidean curvature.

The examples in \cite[\S3]{karpenkov} show that a line has constant $A$, a centred
ellipse has constant $A/B$, and a logarithmic spiral centred at $O$ has
constant $A^3B$. By the curvature identity, $A^2B$ is constant on any
circle.

In \cite[Remark~15.23]{book} Karpenkov notes that it is natural to study the
curves defined by relations between the densities, and asks which curves
are obtained if $AB$, or $B$ alone, is constant. These two questions are
the main subject of this note.

We treat these relations as differential equations for the curve. The
basic observation is that, wherever $A\ne0$ and the real powers are defined, the monomial relation
$A^pB=C$ is equivalent to $\kappa=CA^{2-p}$. This connects the questions
of \cite{book} with several classical subjects, including self-similar solutions
of curve-shortening flow, the support-function equation of convex curves,
the apsidal angle of power-law central forces, and constant centro-affine
curvature. In Section~4 we note that $A\,ds$ and $B\,ds$ are the
centro-affine arclength and turning, so the densities are classical
centro-affine invariants measured per unit of Euclidean length.

The novelty claimed is limited. The densities, the curvature identity and
the questions studied here are due to O.~Karpenkov \cite{karpenkov,book}. Constant
centro-affine curvature, the Abresch--Langer shrinkers, polar reciprocity,
the hodograph, and the monotonicity of the apsidal angle are classical.

The contributions of this note are (i) the sign convention that makes the
curvature identity correct for signed quantities (Proposition~\ref{res:2.1}),
(ii) the identification of $AB=\mathrm{const}$ with the soliton equation of
curve-shortening flow (Proposition~\ref{res:5.1}), (iii) for constant nonzero $B$,
the reduction to the constant-force Binet equation, the interpretation as
polar reciprocals and as hodographs of constant-force orbits, and the
complete closure criterion (Theorems~\ref{res:7.2} and~\ref{res:7.3}, Propositions~\ref{res:7.4}
and~\ref{res:7.5}), and (iv) the explicit local solution of $A^3B=\mathrm{const}$
(Proposition~\ref{res:8.1}). We also describe the ends of the open curves in Sections~5, 6 and~8 and give the turning fractions for general positive support in Section~9. Section~10 treats a general relation $B=\Phi(A)$: the force law (Remark~\ref{res:10.1}), the ends of the arcs with $AB<0$ (Proposition~\ref{res:10.2}), and the number of double points of the closed curves (Proposition~\ref{res:10.3}). It also constructs relations with a prescribed turning fraction (Remark~\ref{res:10.4}), and it describes the curves whose areal density is linear in the arclength and those whose angular density is linear in the arclength (Theorem~\ref{res:10.5} and Corollary~\ref{res:10.6}). Proposition~\ref{res:10.7} gives tests for the monotonicity of the turning fraction of a general relation, among them an explicit test at the circular solution, and Proposition~\ref{res:10.8} shows that for smooth relations no condition on finitely many derivatives, imposed separately at each point of the well, is equivalent to an everywhere positive (or everywhere negative) derivative of the turning fraction.

\section{The densities of O.~Karpenkov and the sign of the curvature identity}

Let $\gamma(s)$ be a $C^2$ regular plane curve parametrized by Euclidean
arclength, and let $O$ be the origin. Write $T=\gamma'$ for the unit
tangent, $N=JT$ for the unit normal, where $J$ is the rotation by
$+\pi/2$, and $\kappa=\det(T,T')$ for the signed curvature, so that
$T'=\kappa N$ and $N'=-\kappa T$. Throughout, $\det$ denotes the signed
area of the parallelogram spanned by two vectors.

\paragraph{Where the densities come from.}
Let $A_0A_1\dots A_n$ be a broken line such that $O$, $A_k$ and $A_{k+1}$
are never collinear. Karpenkov \cite[\S2.2]{karpenkov} assigns to it, with
respect to $O$, a sequence of real numbers, its LLS sequence, in which
elements of two kinds alternate. Each edge $A_kA_{k+1}$ contributes
$\det(A_k,A_{k+1})$, twice the signed area of the triangle $OA_kA_{k+1}$.
Each vertex $A_k$ contributes
$\det(A_{k-1}-A_k,\,A_{k+1}-A_k)$ divided by the contributions of the two
adjacent edges; this measures how sharply the broken line turns at $A_k$,
as seen from $O$. For the sail of a lattice cone these numbers are the
integer lengths of the edges and the integer sines of the angles, and the
continued fraction formed from them is the slope of the cone
\cite[Theorem~1.2]{karpenkov}. Replacing the broken line by a curve and the
elements by their densities per unit length gives the following
definition.

\begin{definition*}[Karpenkov {\cite[Definition~3.1]{karpenkov}},
{\cite[Definition~15.18]{book}}]
The \emph{areal density} and the \emph{angular density} of $\gamma$ with
respect to $O$ at the parameter $s$ are
\[
 A(s)=\lim_{\varepsilon\to0}
 \frac{\det\bigl(\gamma(s),\gamma(s+\varepsilon)\bigr)}{\varepsilon},
 \qquad
 B_K(s)=\lim_{\varepsilon\to0}
 \frac{\det\bigl(\gamma(s-\varepsilon)-\gamma(s),\,
                   \gamma(s+\varepsilon)-\gamma(s)\bigr)}
 {\varepsilon\,\det\bigl(\gamma(s-\varepsilon),\gamma(s)\bigr)\,
                \det\bigl(\gamma(s),\gamma(s+\varepsilon)\bigr)}.
\]
\end{definition*}

Thus, for the broken line with vertices $\gamma(k\varepsilon)$, the edge
elements are $\varepsilon A+o(\varepsilon)$ and, where $A\ne0$, the vertex
elements are $\varepsilon B_K+o(\varepsilon)$. In \cite{karpenkov} the
determinants are written $|v\times w|$, the oriented area of the
parallelogram spanned by $v$ and $w$; in \cite{book} they are written
explicitly. In the numerator of $B_K$ the backward chord is written first;
Proposition~\ref{res:2.1} shows that this order fixes the sign of the
curvature identity.

Since $\gamma(s+\varepsilon)=\gamma(s)+\varepsilon T(s)+o(\varepsilon)$,
the areal density is
\[
 A(s)=\det(\gamma(s),T(s)).
\]
We also put
\[
 X(s)=\det(\gamma(s),N(s))=\langle\gamma(s),T(s)\rangle .
\]
Since $\det(a,b)=\langle Ja,b\rangle$ and $J$ is skew,
$A=-\langle\gamma,N\rangle$. So $|A|$ is the distance from $O$ to the
tangent line at $\gamma(s)$, and $\tfrac12|A|$ is the rate at which the
segment from $O$ to $\gamma(s)$ sweeps out area. For a positively oriented
locally convex arc, $A$ is the Euclidean support function with respect to
$O$.
The point of the curve is recovered as $\gamma=XT-AN$. Differentiation gives
\begin{equation}
 A'=\kappa X,\qquad X'=1-\kappa A. \tag{2.1}
\end{equation}

\begin{proposition}[Sign of the curvature identity]
\label{res:2.1}
At every point where $A\ne0$ the limit $B_K$ exists and
$A^2B_K=-\kappa$. This identity does not depend on the orientation of
the curve.
\end{proposition}

\begin{proof}
By Taylor's formula,
\[
 \gamma(s\pm\varepsilon)-\gamma(s)
 =\pm\varepsilon T+\tfrac12\varepsilon^2\kappa N+o(\varepsilon^2).
\]
Since $\det(T,N)=1$, the numerator equals
\[
 \det(-\varepsilon T,\tfrac12\varepsilon^2\kappa N)
 +\det(\tfrac12\varepsilon^2\kappa N,\varepsilon T)
 +o(\varepsilon^3)=-\varepsilon^3\kappa+o(\varepsilon^3).
\]
Further,
\[
 \det(\gamma(s-\varepsilon),\gamma(s))
 =\det(\gamma(s)-\varepsilon T+o(\varepsilon),\gamma(s))
 =\varepsilon A+o(\varepsilon),
\]
and similarly $\det(\gamma(s),\gamma(s+\varepsilon))=\varepsilon A+
o(\varepsilon)$. Hence the quotient equals
\[
 \frac{-\varepsilon^3\kappa+o(\varepsilon^3)}
      {\varepsilon^3A^2+o(\varepsilon^3)}
 \longrightarrow -\frac{\kappa}{A^2}.
\]
If the orientation is reversed, $\kappa$ and $A$ change sign, so
$-\kappa/A^2$ changes sign. Consequently $B_K$ changes sign and
$A^2B_K=-\kappa$ is preserved. Directly, the two chords in the numerator
are exchanged and each determinant in the denominator changes sign.
\end{proof}

In \cite[Proposition~3.3]{karpenkov} and \cite[Proposition~15.20]{book} the identity is
stated as $A^2B=\kappa$, with $\kappa$ the signed curvature. In \cite{karpenkov} the
sign convention for $\kappa$ is not specified. In \cite{book} the proof uses
$\kappa=x'y''-y'x''$, which is the standard convention
$\kappa=\det(\gamma',\gamma'')$. With this convention and the chord order
of the definition, Proposition~\ref{res:2.1} gives $A^2B_K=-\kappa$. The difference
arises in the last step of the proof in \cite{book}: for the chords
$v_\pm=\pm\varepsilon\gamma'+\tfrac12\varepsilon^2\gamma''$ one has
$\det(v_-,v_+)=-\varepsilon^3(x'y''-y'x'')$, with a minus sign. The
identity $A^2B=\kappa$ holds as stated if the forward chord is written
first, or if curvature is taken positive for clockwise turning.

For example, the counterclockwise unit circle about its centre has
$A=1$, $\kappa=1$ and $B_K=-1$.

The chord order chosen in \cite{karpenkov,book} is natural for continued fractions:
a sail traversed counterclockwise about $O$ turns clockwise, so
$\kappa<0$ and $B_K>0$, in agreement with the positive elements of the
LLS sequence.

The smooth examples in \cite{karpenkov} give positive values of $B$ for curves
traversed counterclockwise. These values are $|B_K|$. For instance, for
the centred ellipse $x^2/a^2+y^2/b^2=1$ traversed counterclockwise one
finds $A/B_K=-a^2b^2$, in agreement with the value $a^2b^2$ of
$|A/B|$ in \cite[Example~3.8]{karpenkov}.

\paragraph{Convention.}
In this note we use the angular density with the opposite chord order,
\begin{equation}
 B:=-B_K,\qquad \text{so that}\qquad A^2B=\kappa. \tag{2.2}
\end{equation}
Every statement below refers to $B$. To translate a statement into the
density $B_K$ of \cite{karpenkov,book}, replace each constant $C$ by $-C$.

Reversing the orientation of the curve changes the signs of $A$, $B$ and
$\kappa$. From (2.2), wherever $A\ne0$,
\begin{equation}
 B=\kappa/A^2,\qquad A^pB=C\ \Longleftrightarrow\ \kappa=CA^{2-p}.
 \tag{2.3}
\end{equation}
We call an arc \emph{nondegenerate} if $A\ne0$ on it, that is, if no
tangent line of the arc passes through $O$. Throughout this paper we
restrict the angular density to nondegenerate arcs, where $B=\kappa/A^2$.

\section{The monomial hierarchy $A^pB=C$}

\begin{center}
\small
\begin{tabular}{@{}c l l p{.42\linewidth}@{}}
\hline
$p$ & Relation & Curvature equation & Interpretation \\
\hline
0 & $B=C$ & $\kappa=CA^2$ & constant angular density (\S\S6--7) \\
1 & $AB=C$ & $\kappa=CA$ & curve-shortening solitons (\S5) \\
2 & $A^2B=C$ & $\kappa=C$ & constant curvature \\
3 & $A^3B=C$ & $\kappa A=C$ & linear equation in the tangent angle (\S8) \\
\hline
\end{tabular}
\end{center}

The case $p=2$ is the curvature identity itself. The cases $p=0$ and
$p=1$ are the two relations asked about in \cite[Remark~15.23]{book}. The case
$p=3$ contains the logarithmic spirals, but its solution space is larger
(\S8).

For completeness we also treat the relation $A=\mathrm{const}\ne0$.
By (2.1), $A'=\kappa X$ and
$X=\langle\gamma,T\rangle=\tfrac12(|\gamma|^2)'$.

Let $U$ be the open subset of a connected nondegenerate arc on which
$\kappa\ne0$. On $U$ one has $X=0$, so $|\gamma|$ is locally constant.
Hence each component of $U$ is an arc of a circle centred at $O$, of
radius $|A|$, and $|\kappa|=1/|A|$ there. By continuity,
$|\kappa|=1/|A|\ne0$ on the closure of $U$, so $U$ is also closed.

Therefore either $U$ is empty, and the arc is a segment of a straight
line at distance $|A|$ from $O$, or $U$ is the whole arc, and the arc
is an arc of the circle of radius $|A|$ centred at $O$.

\begin{example}[Constant curvature and the choice of origin]
The counterclockwise circle $\gamma(t)=(d+R\cos t,R\sin t)$ has
$A=R+d\cos t$ and $\kappa=1/R$. For $0<d<R$ its areal density varies while
$A^2B=1/R$ everywhere. At $d=R$ the support vanishes at $t=\pi$:
the circle stays smooth while the angular density becomes undefined.
\end{example}

\section{The ratio $B/A$ and centred conics}

The ratio $B/A$ serves as a benchmark.

On a nondegenerate arc, let $\sigma$ be the centro-affine arclength,
defined by $d\sigma/ds=A$, so that $\det(\gamma,\gamma_\sigma)=1$ and
$\gamma_\sigma=T/A$. Differentiating $\det(\gamma,\gamma_\sigma)=1$
gives $\det(\gamma,\gamma_{\sigma\sigma})=0$, so
$\gamma_{\sigma\sigma}=-k\gamma$ for a function $k$, the centro-affine
curvature. From
\[
 \gamma_{\sigma\sigma}=A^{-1}(T/A)'
 =A^{-1}\bigl(\kappa N/A-A'T/A^2\bigr)
\]
one obtains
\[
 k=\det(\gamma_\sigma,\gamma_{\sigma\sigma})=\frac{\kappa}{A^3}
   =\frac BA.
\]
Thus constant $B/A$ is exactly constant centro-affine curvature, and
the classification is classical. See O.~Mayer and A.~Myller \cite{mayer}.

\paragraph{The densities as centro-affine invariants.}
The computation above gives a dictionary between Karpenkov's densities
and the centro-affine geometry of curves with respect to $O$:
\[
 A\,ds=\det(\gamma,d\gamma)=d\sigma,\qquad
 B\,ds=\frac BA\,A\,ds=k\,d\sigma .
\]
The two differentials on the right are invariant under every linear map
of determinant~$1$. So are the edge and vertex elements of the LLS
sequence of a broken line, since they are built from determinants
\cite[Proposition~2.3]{karpenkov}; for the broken line with vertices
$\gamma(j\varepsilon)$ they are $\varepsilon A+o(\varepsilon)$ and, up to the
sign fixed in (2.2), $\varepsilon B+o(\varepsilon)$. The densities themselves are
not invariant, because they are taken per unit of Euclidean arclength.
Thus $A$ and $B$ are the centro-affine arclength and the centro-affine
turning $k\,d\sigma$, each measured against Euclidean length, and a relation
$B=\Phi(A)$ combines the two geometries: in centro-affine terms it reads
$k=\Phi(A)/A$ with $A=d\sigma/ds$. What Karpenkov's construction adds to the
classical theory is the discrete origin of these invariants: the
elements of the continued fraction of a broken line become, in the
limit, $d\sigma$ and $k\,d\sigma$.

\begin{remark*}
A different passage from lattice data to an affine invariant of a smooth
convex curve appears in recent work of N.~Kalinin, E.~Lupercio and
M.~Shkolnikov \cite{kls}. They define an $\mathrm{SL}_2(\mathbb Z)$-invariant
tropical zeta function of a convex domain, which in dimension two is a
Dirichlet series over Farey pairs. For $C^3$ strictly convex domains it
has a simple pole at $s=2/3$ whose residue is proportional to the
equiaffine perimeter. There the smooth invariant is a residue of a sum
over all lattice directions, and it does not depend on a centre. Here,
as in \cite{karpenkov,book}, it is a pointwise limit of determinants
taken relative to the fixed point $O$, and the ratio $B/A$ is the
centro-affine curvature. We do not know a formula that relates the two
constructions.
\end{remark*}

\begin{proposition}
\label{res:4.1}
Let $\gamma$ be a nondegenerate arc with $B/A=c$ constant.
If $c>0$, $\gamma$ is an arc of an ellipse centred at $O$.
If $c<0$, $\gamma$ is an arc of one branch of a hyperbola centred at $O$.
If $c=0$, $\gamma$ is an arc of a straight line not passing through $O$.
\end{proposition}

\begin{proof}
The equation is $\gamma_{\sigma\sigma}+c\gamma=0$.
If $c=\omega^2>0$, then
$\gamma=U\cos(\omega\sigma)+V\sin(\omega\sigma)$, and
$\det(\gamma,\gamma_\sigma)=\omega\det(U,V)=1$, so $U$ and $V$ are
linearly independent. The image is a centred ellipse.
If $c=-\lambda^2<0$, then
$\gamma=U\cosh(\lambda\sigma)+V\sinh(\lambda\sigma)$
with $\lambda\det(U,V)=1$, which is one branch of a centred hyperbola.
If $c=0$, then $\gamma=U+\sigma V$ with $\det(U,V)=1$, a straight line
that does not pass through $O$.
\end{proof}

In particular, the centred ellipse of \cite[Example~3.8]{karpenkov} is the case
$c>0$, and for $c\ne0$ the relation can be written $A/B=1/c$.
With the counterclockwise orientation the ellipse $x^2/a^2+y^2/b^2=1$
has $c=1/(a^2b^2)$.

\begin{example}[Constant $B/A$]
The ellipse $(2\cos t,\sin t)$ has $B/A=1/4$ and closes smoothly.
The hyperbola $(\cosh t,\sinh t)$ has $B/A=-1$: its two ends approach
the lines $y=\pm x$. In the latter parametrization
$\det(\gamma,\dot\gamma)=1$ and $\ddot\gamma=\gamma$, so the
centro-affine curvature of Proposition~\ref{res:4.1} is $-1$.
\end{example}

\begin{example}[Centro-affine curvature linear in $\sigma$: the Airy curve]
Let $B/A=c\sigma$ with $c\ne0$ constant, so that
$\gamma_{\sigma\sigma}=-c\,\sigma\gamma$. With $t=-c^{1/3}\sigma$ this is
Airy's equation $\gamma_{tt}=t\gamma$, and
\[
 \gamma=U\operatorname{Ai}(t)+V\operatorname{Bi}(t),\qquad
 -c^{1/3}\det(U,V)=\pi,
\]
because the Wronskian of $\operatorname{Ai}$ and $\operatorname{Bi}$ is
$1/\pi$ and $\det(\gamma,\gamma_\sigma)=1$. The curve is nondegenerate for
all $t$, and it has exactly one inflection point, at $t=0$. As
$t\to+\infty$, $\operatorname{Ai}(t)\to0$ and
$\operatorname{Bi}(t)\to\infty$, so this end goes to infinity and approaches
the ray from $O$ spanned by $V$. As $t=-x\to-\infty$, with
$\zeta=\tfrac23x^{3/2}$ \cite[\S9.7]{dlmf},
\[
 \gamma=\pi^{-1/2}x^{-1/4}\bigl(U\sin(\zeta+\tfrac\pi4)
 +V\cos(\zeta+\tfrac\pi4)\bigr)+O(x^{-7/4}),
\]
so this end winds infinitely often around $O$ along shrinking ellipses and
tends to $O$. Its Euclidean length is infinite, because $|\gamma_t|$ grows
like $x^{1/4}$. In contrast, the logarithmic spiral of Section~8 reaches
$O$ after finite length.
\end{example}

\section{The relation \(AB=\mathrm{const}\) and curve-shortening solitons}
Suppose \(AB=C\) on a nondegenerate arc. By (2.3),
\begin{equation}
\kappa=CA=-C\langle\gamma,N\rangle. \tag{5.1}
\end{equation}
A dilation \(\gamma\mapsto\lambda\gamma\), \(\lambda>0\), maps \(\kappa\) to \(\kappa/\lambda\) and \(A\) to \(\lambda A\), so it replaces \(C\) by \(C/\lambda^2\). The sign of \(C\) is invariant.

A curve \(\gamma_0\) generates a self-similar solution \(\gamma_t=\sqrt{1-2ct}\,\gamma_0\) of curve-shortening flow if and only if \(\kappa=-c\langle\gamma_0,N\rangle\), where \(c>0\) corresponds to shrinking and \(c<0\) to expanding solutions. With the normalization \(c=1/2\), this is the standard shrinker equation \(\kappa=-\tfrac12\langle\gamma,N\rangle\) \cite{al,chou}. Comparing with (5.1) gives the following statement.

\begin{proposition}
\label{res:5.1}
A nondegenerate arc satisfies \(AB=C\) if and only if it satisfies the self-similar equation of curve-shortening flow with \(c=C\). Hence \(C>0\) gives exactly the self-shrinkers, \(C<0\) the self-expanders, and \(C=0\) the straight lines not through \(O\).
\end{proposition}

For \(C>0\) the curvature does not change sign. Differentiating (5.1) with \(C=1\) gives \(\kappa'=\kappa X\), so \(\kappa\) vanishes at one point only if it vanishes identically.
Choose the orientation with \(\kappa>0\) and let \(h(\theta)\) be the support function in the normal-angle parameter \(\theta\). In the normalization \(C=1\) one has \(h=A=\kappa\) and \(1/\kappa=h+h''\), hence
\[
h''+h=\frac1h,\qquad \frac12(h')^2+\frac12h^2-\log h=E .
\]
The closed solutions of this equation are classified by U.~Abresch and J.~Langer \cite{al}. The circle is the only embedded closed shrinker, and the other closed shrinkers are immersed curves with rotation numbers \(m/n\) in the interval \((1/2,1/\sqrt2)\), one for each such rational.

The classification of solitons is classical. The result here is the identification of the relation \(AB=\mathrm{const}\) of \cite{book} with the soliton equation.

\begin{example}[Shrinking circle and stationary line]
A centred circle of initial radius $R$ has $AB=R^{-2}$.
Under curve shortening its radius is $\sqrt{R^2-2t}$.
Its curvature tends to infinity as $t\uparrow R^2/2$, when the circle
collapses to $O$. For the oriented line $\gamma(s)=(s,d)$, $d\ne0$,
we have $A=-d$ and $B=0$. The line is stationary.
\end{example}
\begin{example}[An expanding graph with two asymptotic rays]
Take $C=-1$ and write the profile as $(x,f(x))$, oriented to the right.
Equation (5.1) becomes
\[
 f''=(1+f'^2)(f-xf'),\qquad f(0)=1/2,\quad f'(0)=0.
\]
The solution is even by uniqueness. Put $g=f-xf'$. For $x\ge0$,
$g'=-x(1+f'^2)g$, so $0<g\le \tfrac12 e^{-x^2/2}$.
Since $(\arctan f')'=g$, the tangent angle stays below
$\tfrac12\sqrt{\pi/2}<\pi/2$. The slope stays bounded, giving a global
strictly convex graph. Its limiting positive slope $m$ exists;
$f''=O(e^{-x^2/2})$ implies $x(m-f'(x))\to0$, and hence
$f(x)-mx=g-x(m-f'(x))\to0$. By evenness the asymptotic rays are
$y=m|x|$. The support $A=-g/\sqrt{1+f'^2}$ never vanishes at a finite
point. The flow dilates this smooth profile by $\sqrt{1+2t}$ for
$t>-1/2$. In the backward limit the two asymptotic rays meet at $O$.
\end{example}

\paragraph{The full profiles and the time evolution.}
On every nonstraight profile of (5.1),
$A_s=\kappa\langle\gamma,T\rangle$ and $\kappa=CA$ give
\[
 |\kappa(s)|=K\exp\!\big(C|\gamma(s)|^2/2\big),\qquad K>0 .
\]
In particular a nonstraight solution has no inflection point.
For $C>0$ the maximal profiles are smooth and bounded. The noncircular ones
consist of successive rotations of one oscillation; the image is closed
when the fraction of a turn per oscillation is rational and dense in an
annulus when it is irrational. Every shrinking profile scales by
$\sqrt{1-2Ct}$, so curvature blows up at extinction even when the initial
profile is smooth. For $C<0$ every maximal nonstraight profile is a smooth
embedded graph after rotation, with two ends asymptotic to rays from $O$.
The curvature tends to zero at infinity, although $|B|=|C/A|$ tends to
infinity there. For the global classification see Halldorsson \cite[Theorems 5.1 and 6.1]{halldorsson}.

\section{Constant angular density \texorpdfstring{$B=\mathrm{const}$}{B = const}}

Suppose that $B=C\ne0$. Reversal of orientation changes the sign of $B$,
so choose the orientation for which $C>0$. On every nondegenerate arc,
$\kappa=CA^2>0$. The curve is therefore locally strictly convex and its
normal angle $\theta$ is a regular parameter. Write
$e(\theta)=(\cos\theta,\sin\theta)$ for the outward unit normal.
The signed support function $h=A$ satisfies $1/\kappa=h+h''$, so
\[
  h''+h=\frac{1}{Ch^2}.
\]
This equation holds for either sign of $h$. We first consider $h>0$;
Lemma~\ref{res:6.1} below shows that this includes every closed nondegenerate curve.
Writing $h=C^{-1/3}u$ gives
\begin{equation}
  u''+u=u^{-2}, \qquad u>0. \tag{6.1}
\end{equation}
Multiplication by $u'$ gives the conserved energy
\begin{equation}
  \tfrac12 (u')^2+\tfrac12 u^2+\frac1u=E,\qquad E\ge\tfrac32.
  \tag{6.2}
\end{equation}
For nonconstant positive periodic solutions let $P(E)$ denote the least period in normal angle.

\begin{lemma}
\label{res:6.1}
Let $\gamma$ be a closed regular $C^2$ curve on which $A$ never vanishes,
and suppose that $B=C>0$. Then $\gamma$ is locally strictly convex,
and its signed support function $h=A$ is positive.
\end{lemma}

\begin{proof}
We have $\kappa=CA^2>0$. Since $A$ is continuous and nonzero on the
connected curve, it has constant sign. During one traversal the tangent
angle, and hence the normal angle, increases by $2\pi m$ for some integer
$m\ge1$. Choose the initial angle as the lower integration limit
$\theta_0$ and put $\theta_1=\theta_0+2\pi m$. The derivative $h'$ returns
to its initial value, so
\[
 \int_{\theta_0}^{\theta_1} h\,d\theta
 =\int_{\theta_0}^{\theta_1}(h+h'')\,d\theta
 =\int_{\theta_0}^{\theta_1}\frac{d\theta}{\kappa}
 =\operatorname{length}(\gamma)>0.
\]
Thus $h=A>0$ throughout the curve.
\end{proof}

\begin{proposition}[Reconstruction and regularity for constant angular density]
\label{res:6.2}
Fix $C>0$. On the branch $h=A>0$, curves with $B=C$ correspond locally,
in outward normal angle $\theta$, to positive solutions of
$u''+u=u^{-2}$ by $h=C^{-1/3}u$ and
$\gamma=h e+h'e'$, where $e=(\cos\theta,\sin\theta)$.
Every maximal such solution extends to all real $\theta$. The equilibrium
is the centred circle of radius $C^{-1/3}$. At every energy $E>3/2$
the support oscillates periodically between two positive turning points,
and the reconstructed curve is a regular immersion.
\end{proposition}
\begin{proof}
The support identity gives $h+h''=1/(Ch^2)$; conversely
$\gamma'=(h+h'')e'=C^{-1/3}u^{-2}e'$ never vanishes and gives
$\kappa=Ch^2$. The potential $V(u)=u^2/2+1/u$ is strictly convex on
$(0,\infty)$, tends to infinity at both ends, and has minimum $3/2$ at
$u=1$. The conserved energy bounds $u$ away from zero and infinity and
bounds $u'$, giving global existence. For $E>3/2$ both turning points
are simple, so the time integral is finite at each and the motion repeats.
\end{proof}
\begin{example}[Two turning points do not imply closure]
For $C=1$, take $u(0)=2$, $u'(0)=0$. Then $E=5/2$ and the turning
points are $2$ and $\sqrt2-1$, the positive roots of $u^3-5u+2=0$.
The curvature $u^2$ is finite and positive throughout. Whether the
plane curve closes depends on the ratio of the support period to $2\pi$,
as in Section~7. Figure~\ref{fig:closedfour} there supplies a closed example.
\end{example}

\begin{proposition}[The two signs of the support for $B=C$]
\label{res:6.3}
Orient the curve so that $C>0$. Use the right unit normal
$e=(\cos\theta,\sin\theta)$, so that $e'=T$ and $h=A$.
For positive support, a nonconstant solution with irrational $P/(2\pi)$
has image dense in a closed annulus. For negative support, every maximal
solution is a smooth embedded curve of infinite length in each direction.
Its two ends approach rays from $O$, and its curvature tends to zero
along each end.
\end{proposition}
\begin{proof}
Work with $h$ itself, before rescaling. Multiplying the support equation
by $h'$ gives
\[
 E=\tfrac12(h'^2+h^2)+\frac1{Ch},\qquad
 |\gamma|^2=h^2+h'^2=2E-\frac2{Ch}.
\]
If $h>0$ oscillates between $a$ and $b$, then $|\gamma|$ does too:
its square increases with $h$ and equals $h^2$ at the turning points.
After each support period the whole arc rotates by $P$. If $P/(2\pi)$
is irrational, these rotations are dense on every circle whose radius
is attained during one oscillation. Their union is the annulus
$a\le |x|\le b$. The speed $1/(Ch^2)$ in normal angle and curvature
$Ch^2$ remain positive and bounded.

If $h<0$, the potential $h^2/2+1/(Ch)$ decreases strictly from
$+\infty$ to $-\infty$. There is one turning point, a minimum.
On either side $h$ tends to $0$ from below and
$|h'|\sim\sqrt{2/(C|h|)}$. The angular distance to the endpoint
is therefore $O(|h|^{3/2})$. The radius $|\gamma|$ tends to infinity,
so each end has infinite length. If $e_*$ is the endpoint normal,
$\gamma\cdot e_*=O(|h|)+O(|h'|\,|h|^{3/2})\to0$.
Thus each end approaches a ray of the line through $O$ perpendicular
to $e_*$. Its curvature $Ch^2$ tends to zero.
To see embeddedness, put the minimum at $\theta=0$ and let $g=-h$.
Then $g'(0)=0$ and $g''+g=-1/(Cg^2)<0$.
Variation of constants gives
$g(t)=g(0)\cos t-\int_0^t\sin(t-v)/(Cg(v)^2)\,dv$.
Hence the positive interval ends before $\pi/2$ in either direction
(by evenness). The tangent turns through less than $\pi$ in total,
and its projection onto the middle tangent direction stays positive.
\end{proof}
\begin{example}
Set $C=1$, $h(0)=-1$, $h'(0)=0$. Then $E=-1/2$ and the preceding
identity gives $h=-2/(1+|\gamma|^2)$ and
$\kappa=4/(1+|\gamma|^2)^2$. The closest point has distance one from $O$.
\end{example}

\section{The Binet equation, the period bound, and closure}
Equation (6.1) is the orbit equation of a central force, known as Binet's equation, for a particle of unit mass moving under an attractive central force of constant magnitude \(F\), with angular momentum \(L\), in the normalization \(F=L^2=1\). If \(r\) is the distance to the centre and \(w=1/r\), Binet's equation reads \(w''+w=F/(L^2w^2)\).

In this normalization the mechanical energy is \(\tfrac12((w')^2+w^2)+1/w\), so \(E\) in (6.2) is exactly the energy of the orbit. Half of the least period \(P(E)\) is the apsidal angle \(\Theta(E)\) of the orbit. In this first-integral form the statement is in V.~I.~Arnold \cite[\S8, Problem~1]{arnold}: for a potential \(U(r)\) and angular momentum \(M\), the angle between a pericentre and the next apocentre is half the period of the one-dimensional system with potential energy \(W(x)=U(M/x)+x^2/2\), where \(x=M/r\). For \(U=r\) and \(M=1\) this is (6.2).

In the notation of R.~Castelli \cite{castelli} and D.~Rojas \cite{rojas}, the force is \(f(r)\sim r^{-(\alpha+1)}\) and the potential is \(V(r)=-\alpha^{-1}r^{-\alpha}\). The constant force corresponds to \(\alpha=-1\).
Castelli proved that for \(-2<\alpha<1\) the apsidal angle is a monotone function of the eccentricity. Rojas \cite[Theorem~1.1]{rojas} proved that, for fixed nonzero angular momentum, \(\Theta'(E)<0\) for \(\alpha\in(-2,-2/3)\).
Since \(\alpha=-1\) lies in this interval, \(P(E)\) is a strictly decreasing function of \(E\) on \((3/2,\infty)\). It remains to compute the two limits.

Near the circular orbit, put \(u=1+v\). Since \(u^{-2}=1-2v+O(v^2)\), equation (6.1) becomes \(v''+3v=O(v^2)\). Hence \(P(E)\to2\pi/\sqrt3\) as \(E\downarrow3/2\).
For large energies, by (6.2),
\[
P(E)=2\int_{u_-}^{u_+}\left(2E-u^2-\frac2u\right)^{-1/2}\,du .
\]
Substitute \(u=\sqrt{2E}\,x\) and put \(\varepsilon=2(2E)^{-3/2}\). Then
\[
P(E)=2\int_{x_-}^{x_+} f_\varepsilon(x)^{-1/2}\,dx,\qquad
f_\varepsilon(x)=1-x^2-\frac{\varepsilon}{x},
\]
where \(x_-<x_+\) are the zeros of \(f_\varepsilon\) in \((0,1)\). As \(\varepsilon\to0\), \(x_-=\varepsilon+O(\varepsilon^3)\) and \(x_+\to1\). Fix \(\delta\in(0,1/4)\) and split the integral at \(2\varepsilon\) and \(1-\delta\).

\noindent (a) On \([x_-,2\varepsilon]\), put \(x=\varepsilon y\). Then \(f_\varepsilon(\varepsilon y)=1-1/y-\varepsilon^2y^2\), and the integral equals
\[
\varepsilon\int_{y_-}^{2}(1-1/y-\varepsilon^2y^2)^{-1/2}\,dy
\]
with \(y_-\to1\). The integrand has an integrable square-root singularity at \(y_-\) and is bounded elsewhere, so this part is \(O(\varepsilon)\).

\noindent (b) On \([2\varepsilon,1-\delta]\) one has \(\varepsilon/x\le1/2\). Choose \(x_0>0\) with \(x_0^2<1/4\). For \(x\le x_0\), \(f_\varepsilon(x)\ge1/2-x_0^2>1/4\). For \(x_0\le x\le1-\delta\), \(f_\varepsilon(x)\ge1-(1-\delta)^2-\varepsilon/x_0\), which is bounded below by a positive constant for small \(\varepsilon\). Thus the integrand is bounded uniformly in \(\varepsilon\) and converges pointwise to \((1-x^2)^{-1/2}\). By dominated convergence, this part tends to \(\int_0^{1-\delta}(1-x^2)^{-1/2}\,dx\).

\noindent (c) On \([1-\delta,x_+]\), \(f_\varepsilon'(x)=-2x+\varepsilon/x^2\le-1/2\) for small \(\varepsilon\), and \(f_\varepsilon(x_+)=0\). Hence \(f_\varepsilon(x)\ge(x_+-x)/2\), and this part is at most \(\int_0^\delta\sqrt2\,t^{-1/2}\,dt=2\sqrt{2\delta}\), uniformly in \(\varepsilon\).

Letting first \(E\to\infty\) (that is, \(\varepsilon\to0\)) and then \(\delta\to0\), we obtain \(P(E)\to2\int_0^1(1-x^2)^{-1/2}\,dx=\pi\). Both limits agree with \cite[\S8, Problems~2--4]{arnold}: for the potential \(U=ar^{\beta}\) the apsidal angle of a nearly circular orbit is \(\pi/\sqrt{\beta+2}\), which is \(\pi/\sqrt3\) for \(\beta=1\), and for \(\beta>0\) it tends to \(\pi/2\) as \(E\to\infty\). Together with strict monotonicity this proves the period bound.

\begin{proposition}[Period bound]
\label{res:7.1}
The least period \(P(E)\) of a nonconstant solution of (6.1) is a continuous, strictly decreasing function of \(E\in(3/2,\infty)\), with \(P(E)\to2\pi/\sqrt3\) as \(E\downarrow3/2\) and \(P(E)\to\pi\) as \(E\to\infty\). In particular,
\begin{equation}
\pi<P(E)<\frac{2\pi}{\sqrt3},\qquad
\text{equivalently}\qquad \frac{\pi}{2}<\Theta(E)<\frac{\pi}{\sqrt3}. \tag{7.1}
\end{equation}
\end{proposition}

Continuity of \(P\) follows from the integral formula above, since \(u_\pm\) depend smoothly on \(E\) for \(E>3/2\).

\begin{theorem}[Embedded curves of constant angular density]
\label{res:7.2}
Let \(\gamma\) be a closed embedded \(C^2\) plane curve on which \(A\) never vanishes. If the angular density \(B\) of \(\gamma\) is a nonzero constant, then \(\gamma\) is a circle centred at \(O\).
\end{theorem}
\begin{proof}
Reverse the orientation if necessary, so that \(B=C>0\). By Lemma~\ref{res:6.1}, \(\kappa>0\) and \(h=A>0\). By the Umlaufsatz of H.~Hopf \cite{hopf}, an embedded closed regular curve has rotation index \(\pm1\). Since \(\kappa>0\), the tangent angle increases by exactly \(2\pi\), so \(\gamma\) is convex and the normal angle \(\theta\) is a global parameter on \([0,2\pi)\). The scaled support function \(u=C^{1/3}h\) is therefore a positive \(2\pi\)-periodic solution of (6.1). Suppose \(u\) is not constant, with least period \(P\). Then \(2\pi=nP\) for a positive integer \(n\). By (7.1), \(n=1\) would give \(P=2\pi>2\pi/\sqrt3\), \(n=2\) would give \(P=\pi\), and \(n\ge3\) would give \(P\le2\pi/3<\pi\). All three are impossible. Hence \(u\) is constant, \(u\equiv1\) by (6.1), and \(\gamma\) is the circle of radius \(C^{-1/3}\) centred at \(O\).
\end{proof}

\begin{remark*}
For a closed convex curve, the condition \(A\ne0\) everywhere is equivalent to \(O\) lying in the interior of the region bounded by \(\gamma\), because a point on the curve or outside it lies on some tangent line.
\end{remark*}

\begin{theorem}[Closed immersed curves of constant angular density]
\label{res:7.3}
Let \(m,n\) be coprime positive integers. A closed immersed nondegenerate curve with \(B=1\) and rotation index \(m\) that is not a circle exists if and only if \(1/2<m/n<1/\sqrt3\), where \(n\) is the number of oscillations of its support function. For each such pair it is unique up to rotation about \(O\), and it corresponds to the unique energy \(E\) with \(P(E)=2\pi m/n\).
\end{theorem}
\begin{proof}
Let \(\gamma\) be such a curve. By Lemma~\ref{res:6.1} its support function \(u\) is a positive solution of (6.1), defined for \(\theta\in\mathbb R\) and \(2\pi m\)-periodic. If \(u\) is not constant, its least period \(P\) satisfies \(2\pi m=nP\) for the number \(n\) of oscillations in one traversal. Since \(\gamma\) does not close at an earlier value \(\theta=2\pi k\) with \(0<k<m\), the integers \(m\) and \(n\) are coprime. By (7.1), \(m/n=P/(2\pi)\in(1/2,1/\sqrt3)\). Conversely, if \(m/n\) lies in this interval, Proposition~\ref{res:7.1} gives a unique \(E\) with \(P(E)=2\pi m/n\). The curve \(\gamma(\theta)=u(\theta)e(\theta)+u'(\theta)e'(\theta)\) satisfies \(B=1\) (\S6) and closes for the first time at \(\theta=2\pi m\), because \(2\pi k\) is a multiple of \(P\) only when \(nk/m\) is an integer, that is, when \(m\) divides \(k\). Uniqueness up to rotation follows from the uniqueness of \(E\) and the translation invariance of (6.1) in \(\theta\).
\end{proof}

\begin{figure}[htbp]
\centering
\includegraphics[width=.55\linewidth]{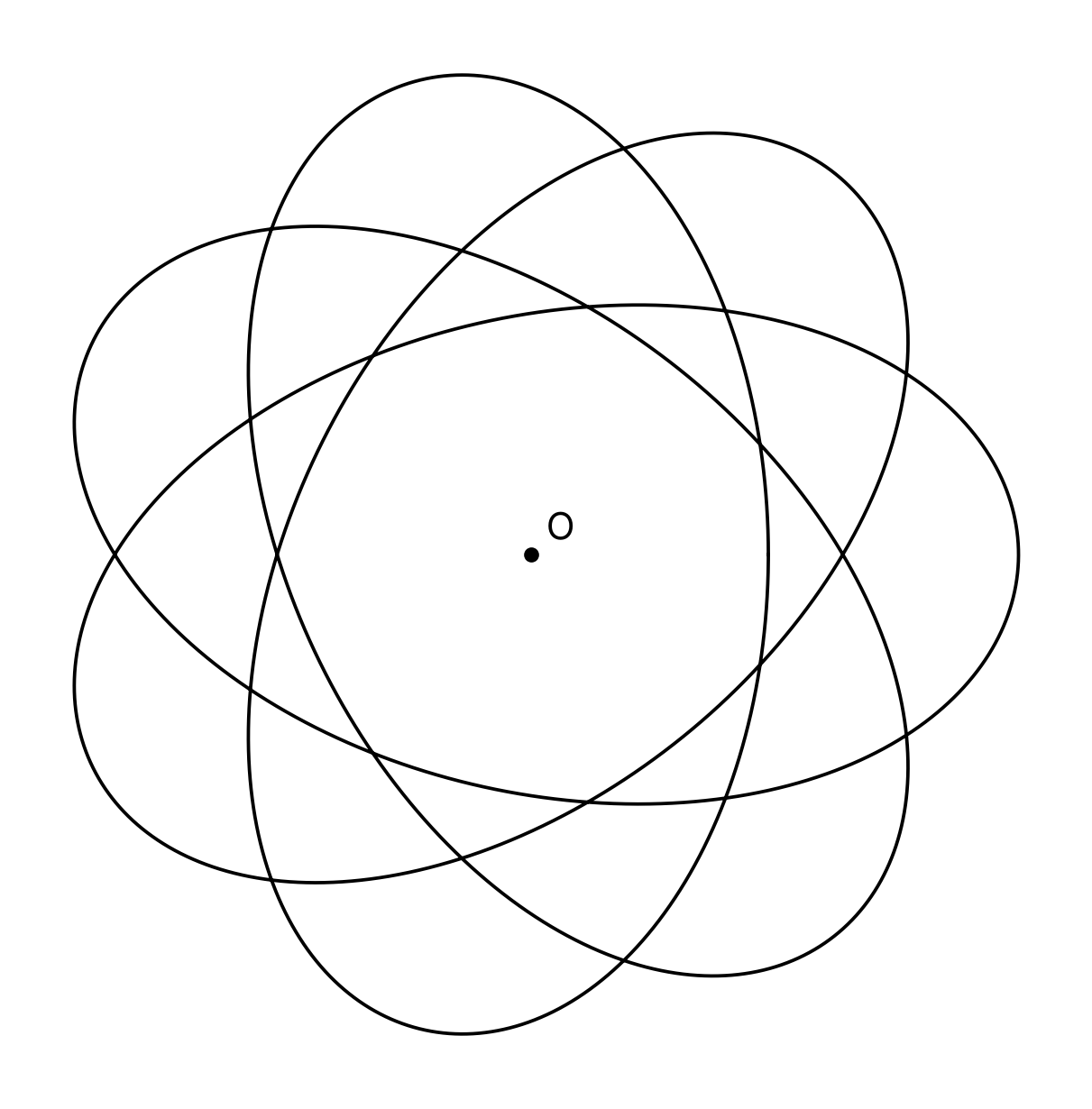}
\caption{The closed curve with constant angular density \(B=1\) and rotation number \(m/n=4/7\) (\(E=1.697378876\ldots\)). The point marks the origin \(O\).}
\label{fig:closedfour}
\end{figure}

The rationals in \((1/2,1/\sqrt3)\) with the smallest numerators are \(4/7\), \(5/9\) and \(6/11\). For \(m/n=4/7\) the energy is \(E=1.697378876\ldots\) (Figure~\ref{fig:closedfour}). By Proposition~\ref{res:10.3} the curve with turning fraction \(m/n\) has exactly \(n(m-1)\) double points; the curve in Figure~\ref{fig:closedfour} has \(21\).
This family is the analogue of the Abresch--Langer family of \S5, whose rotation numbers fill the interval \((1/2,1/\sqrt2)\). There the linearization at the circle is \(v''+2v=0\) instead of \(v''+3v=0\).

\paragraph{Numerical method.} Equation (6.1) was integrated from the inner turning point \(u_-\), with \(u'=0\), by an explicit Runge--Kutta method of order 8 (Dormand--Prince) with relative and absolute tolerances \(10^{-13}\). The period \(P(E)\) is the time of the first return to \(u_-\).
The energy with \(P(E)=8\pi/7\) was found by Brent's method, giving \(E=1.697378876138\). The reconstructed curve closes after \(\theta=8\pi\) with error \(1.1\cdot10^{-12}\), and its angular density, computed by finite differences, equals 1 within \(5\cdot10^{-6}\).
On a grid of 1000 energies in \((3/2,200]\) the computed \(P(E)\) is strictly decreasing, from \(2\pi/\sqrt3=3.627599\) near \(E=3/2\) to \(P(200)=3.14369\). Also \(P(3000)=3.14165\). These values agree with Proposition~\ref{res:7.1}.

The constant-force problem also gives two geometric descriptions of the curves with constant \(B\).
\begin{proposition}[Polar reciprocals of constant-force orbits]
\label{res:7.4}
Let \(\gamma\) be a locally strictly convex arc whose support function \(h(\theta)\) with respect to \(O\) is positive, \(\theta\) being the normal angle. The polar reciprocal \(\gamma^*\) of \(\gamma\) with respect to the unit circle centred at \(O\) has the polar equation \(r(\theta)=1/h(\theta)\). Consequently \(\gamma\) has constant angular density \(B=C>0\) if and only if \(\gamma^*\) is an arc of an orbit of an attractive central force of constant magnitude \(F\) about \(O\) with \(F/L^2=1/C\), where \(L\) is the angular momentum per unit mass.
\end{proposition}
\begin{proof}
The tangent line of \(\gamma\) with normal angle \(\theta\) is
\[
\{x:\langle x,e(\theta)\rangle=h(\theta)\}
=\{x:\langle x,e(\theta)/h(\theta)\rangle=1\}.
\]
Its pole with respect to the unit circle is the point \(e(\theta)/h(\theta)\), and these points form \(\gamma^*\). Put \(w=1/r=h\). For an attractive central force of constant magnitude \(F\), Binet's equation is \(w''+w=(F/L^2)w^{-2}\). This coincides with \(h''+h=1/(Ch^2)\) exactly when \(F/L^2=1/C\).
\end{proof}

\begin{proposition}[Hodographs of constant-force orbits]
\label{res:7.5}
Let \(x(t)\) be a motion of a unit mass under an attractive central force of constant magnitude \(F\) about \(O\), with angular momentum \(L>0\), and let \(v(t)=\dot x(t)\). Then the hodograph \(v\), traversed in the direction of increasing \(t\), is locally strictly convex with positive support function, and its angular density with respect to \(O\) is the constant \(B=1/(FL)\). Conversely, every arc with constant angular density \(B=C>0\) and positive support function is, after a rotation about \(O\), an arc of the hodograph of such a motion with \(F=1/C\) and \(L=1\).
\end{proposition}
\begin{proof}
Write \(x=r e(\varphi)\), so that \(r^2\dot\varphi=L\), \(\dot r=\langle v,e(\varphi)\rangle\) and \(r\dot\varphi=\langle v,Je(\varphi)\rangle=L/r\). Since \(\dot v=-F e(\varphi)\), we have
\[
\frac{dv}{d\varphi}=\frac{\dot v}{\dot\varphi}=-\frac{Fr^2}{L}e(\varphi).
\]
Hence the unit tangent of the hodograph is \(-e(\varphi)\). Its tangent angle \(\varphi+\pi\) increases with \(t\), so its curvature is positive, and its radius of curvature with respect to the tangent angle is \(|dv/d\varphi|=Fr^2/L\). Its unit normal is \(N=J(-e(\varphi))=-Je(\varphi)\), so its areal density is \(A=-\langle v,N\rangle=\langle v,Je(\varphi)\rangle=L/r>0\). Therefore \(\kappa=L/(Fr^2)=A^2/(FL)\), and \(B=\kappa/A^2=1/(FL)\).

For the converse, let \(\gamma\) have \(B=C>0\) and support function \(h>0\) in the normal angle \(\theta\). By Proposition~\ref{res:7.4} with \(L=1\) and \(F=1/C\), the curve \(r=1/h(\varphi)\) is an orbit. By the first part, the support function of its hodograph, as a function of the normal angle \(\varphi+\pi/2\), is \(L/r=h(\varphi)\). A convex arc is determined by its support function, so the hodograph is the image of \(\gamma\) under the rotation by \(\pi/2\).
\end{proof}

Propositions~\ref{res:7.4} and \ref{res:7.5} express the classical fact that the hodograph of a central-force motion is the polar reciprocal of the orbit, rotated by a right angle and scaled by \(L\).
Since \(B\) is invariant under rotations about \(O\), the rotation plays no role for the density. In the Kepler problem the hodograph is a circle. For the constant force it is a curve of constant angular density.
Thus the closed curves with constant \(B\) are exactly the hodographs of closed constant-force orbits, and Theorems~\ref{res:7.2} and \ref{res:7.3} are reformulations of the apsidal-angle analysis.

\section{The relation \texorpdfstring{$A^3B=\mathrm{const}$}{A cubed B = const}}

By (2.3), $A^3B=C$ is equivalent to $\kappa A=C$. For $C=0$,
nondegeneracy gives $\kappa=0$, so the solutions are straight lines
whose supporting line does not pass through $O$.

Assume $C\ne0$. Then $\kappa\ne0$ and the tangent angle $\theta$,
with $d\theta/ds=\kappa$, is a regular parameter. Reversing orientation
changes the signs of $A$ and $\kappa$ together and preserves $C$.
Choose the orientation with $\kappa>0$. From (2.1) and $\kappa=C/A$,
\begin{equation}
 A_\theta=X,\qquad
 X_\theta=\frac{1-\kappa A}{\kappa}=\frac{1-C}{C}A,\qquad
 A_{\theta\theta}=\frac{1-C}{C}A. \tag{8.1}
\end{equation}
In the rest of this section, after Proposition~\ref{res:8.1}, a prime
denotes the derivative with respect to $\theta$, not $s$: $A'=A_\theta$,
$\gamma'=\gamma_\theta$.

\begin{proposition}
\label{res:8.1}
Fix $C\ne0$. For a nondegenerate arc oriented with $\kappa>0$,
the relation $A^3B=C$ holds if and only if the areal density, as a
function of tangent angle, satisfies (8.1) and $A/C>0$.
Conversely, every solution of the last equation in (8.1) on an interval
where $A/C>0$ determines such an arc by
\[
 \gamma=A_\theta T(\theta)-A N(\theta),\qquad
 T(\theta)=(\cos\theta,\sin\theta),\quad N(\theta)=JT(\theta).
\]
The solutions have the following forms, with $\alpha,\beta\in\mathbb R$.
\begin{enumerate}
\item[(i)] If $0<C<1$, then $A=\alpha e^{q\theta}+\beta e^{-q\theta}$,
where $q=\sqrt{(1-C)/C}$.
\item[(ii)] If $C=1$, then $A=\alpha+\beta\theta$.
\item[(iii)] If $C>1$ or $C<0$, then
$A=\alpha\cos(q\theta)+\beta\sin(q\theta)$,
where $q=\sqrt{(C-1)/C}$.
\end{enumerate}
\end{proposition}

\begin{proof}
Necessity follows from the calculation above.
The sign condition $A/C>0$ is equivalent to $\kappa=C/A>0$.
Conversely, let $A$ solve (8.1) with $A/C>0$ and define
$\gamma=A_\theta T-A N$. Since $T_\theta=N$ and $N_\theta=-T$,
\[
 \gamma_\theta=A_{\theta\theta}T+A_\theta N-A_\theta N+AT
 =(A_{\theta\theta}+A)T=(A/C)T.
\]
Thus $T$ is the unit tangent, $ds/d\theta=A/C$, and
$\kappa=d\theta/ds=C/A$. Moreover,
$\det(\gamma,T)=\det(-AN,T)=A$. Hence $A$ is the areal density
of the reconstructed arc and $\kappa A=C$.
The three displayed families are the general solutions of the
constant-coefficient equation for $A$ in (8.1).
\end{proof}

For $0<C<1$, the solutions with $\alpha\beta=0$ are logarithmic spirals
centred at $O$. For the spiral $r=ae^{b\varphi}$ one obtains
$C=1/(1+b^2)$, in agreement with \cite[Example~3.10]{karpenkov}.
The solutions with $\alpha\beta\ne0$ give further curves in this family.

For $C=1$ and $\beta\ne0$, the radius of curvature $A/C=\alpha+\beta\theta$
is linear in the tangent angle, the intrinsic equation of an involute
of a circle. The case $\beta=0$ gives a circle centred at $O$, as in
Section~3. In case (iii), $A$ has zeros and the nondegenerate solutions
are arcs on intervals between consecutive zeros with $A/C>0$.

\begin{example}[A cusp at a zero of the support]
For $C=1$ take $A(\theta)=\theta$ on $\theta>0$ in (8.1).
With $T=(\cos\theta,\sin\theta)$ and $N=(-\sin\theta,\cos\theta)$,
the reconstructed involute is $\gamma=T-\theta N$ and
$\gamma'=\theta T$. Near zero its coordinates are
$(1+\theta^2/2+O(\theta^4),\ \theta^3/3+O(\theta^5))$.
The analytic continuation across zero has an ordinary semicubical cusp;
on the original arc the curvature $1/\theta$ diverges at its endpoint.
\end{example}

\begin{proposition}[Ends of the maximal arcs with $A^3B=C$]
\label{res:8.2}
Use the tangent angle $\theta$ and orient the arc so $\kappa>0$.
Write $T=(\cos\theta,\sin\theta)$, $N=(-\sin\theta,\cos\theta)$.
For $C\ne0$, choose the connected interval on which $A/C>0$.
Then
\[
 \gamma=A'T-AN,\qquad \gamma'=(A/C)T,\qquad
 \kappa=C/A,\qquad |\gamma|^2=A^2+A'^2 .
\]
\begin{itemize}
\item For $0<C<1$ put $q=\sqrt{(1-C)/C}$. Up to translating or
reversing the angle variable, the possible positive functions are
$H\cosh(q\theta)$, $H\sinh(q\theta)$ on $\theta>0$, and
$H e^{q\theta}$, with $H>0$. The first has two smooth ends at infinity.
The second has a finite cusp endpoint and one end at infinity.
The third is a logarithmic spiral: its end at $O$ has finite length
and infinitely many tangent turns; its other end extends to infinity.
\item For $C=1$, constant positive $A$ gives a centred circle.
Nonconstant $A=a+b\theta$ on $A>0$ gives an involute of a circle
with one cusp endpoint and one unbounded end.
\item For $C>1$ or $C<0$, $A$ is sinusoidal. Each maximal arc lies
between two consecutive simple zeros and has finite length and two
cusp endpoints (which need not be distinct in the plane).
\item For $C=0$, all nondegenerate maximal profiles are lines off $O$.
\end{itemize}
\end{proposition}
\begin{proof}
The formulas for $A$ follow from $A''=(1-C)A/C$. At a finite zero
$\theta_0$, $A'(\theta_0)\ne0$ by uniqueness. On extending the formula for
$\gamma$ past the endpoint, $\gamma'(\theta_0)=0$,
$\gamma''=(A'/C)T$ there, and
$\det(\gamma'',\gamma''')=2(A'/C)^2\ne0$ there. Thus in tangent and
normal coordinates the first nonzero terms have orders two and three,
the ordinary cusp. On its regular side the curvature diverges.
For the exponential solution, the length from $-\infty$ to any finite
angle is $\int_{-\infty}^{\theta}He^{qv}/C\,dv<\infty$, although
the tangent angle has no finite limit. The other endpoint assertions
follow from $|\gamma|^2=A^2+A'^2$ and $ds/d\theta=A/C$.
\end{proof}

\begin{example}[Three open curves with the same value of $A^3B$]
For $C=1/2$, equation (8.1) is $A''=A$.
Taking $A=\cosh\theta$ gives a smooth curve with two unbounded ends.
It crosses itself: its first coordinate is odd, its second is even, and
the first vanishes at some $\theta\in(\pi/2,\pi)$ because it changes sign
there. The two tangents at $\theta$ and $-\theta$ are distinct.
Taking $A=\sinh\theta$ on $\theta>0$ gives a finite endpoint at $(1,0)$.
Relative to that point the first terms of the position are
$(\theta^2,\tfrac23\theta^3)$, giving a cusp on continuation.
Taking $A=e^\theta$ gives an end at $O$ of length
$\int_{-\infty}^0 2e^\theta\,d\theta=2$ up to $\theta=0$.
The tangent makes infinitely many turns along this finite length.
\end{example}
\begin{figure}[htbp]
\centering
\includegraphics[width=\linewidth,trim=0 0 0 122,clip]{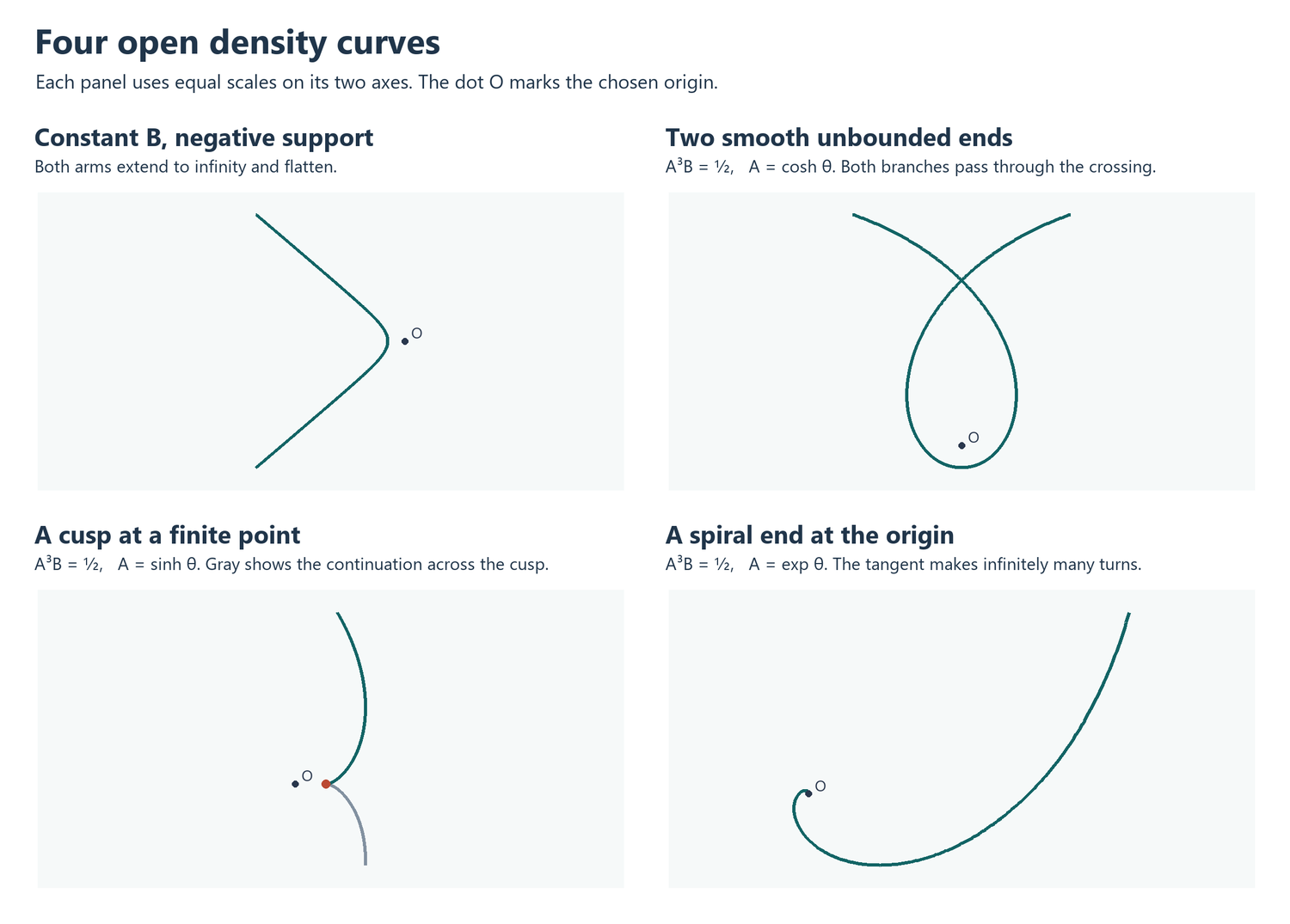}
\caption{Examples from Sections 6 and 8. Each panel uses equal scales on its two axes, and the dot marks the origin $O$. The upper left curve has $B=1$ and negative support. The other three have $A^3B=1/2$. In the lower left panel the gray continuation shows the cusp at the endpoint of the positive-support arc. Only finite portions of the unbounded curves are shown.}
\label{fig:open-curves}
\end{figure}

\section{The family $A^pB=C$}

\paragraph{Physical interpretation.}
On a positive-support arc with $C>0$, the relation $A^pB=C$ gives
$h''+h=C^{-1}h^{p-2}$. Setting $r=1/h$ identifies this with Binet's equation
for an attractive central force of magnitude $F(r)=(L^2/C)r^{-p}$,
where $L$ is the angular momentum per unit mass. The polar reciprocal
of the density curve is therefore a central-force orbit. For nonconstant
periodic support, half its least period is the apsidal angle of that orbit.

The cases $p=0,1,2,3$ correspond respectively to a constant force,
a force proportional to $1/r$, the inverse-square force, and the inverse-cube
force. For each exponent the support of the velocity curve is $L/r$,
as in Proposition~\ref{res:7.5}; after a rotation and a change of scale it recovers
the original density curve.

\paragraph{A test for regularity on positive-support arcs.}
For $C>0$, $h=A>0$ and an arbitrary real $p$, the speed in normal
angle and the curvature are
\[
 \frac{ds}{d\theta}=\frac{h^{p-2}}{C},\qquad \kappa=C h^{2-p}.
\]
Both are positive and finite at every interior point with finite positive
support. The first integral is
\[
 E=\tfrac12 h'^2+V_p(h),\qquad
 V_p(h)=
 \begin{cases}
 \tfrac12h^2-h^{p-1}/[C(p-1)],&p\ne1,\\
 \tfrac12h^2-C^{-1}\log h,&p=1 .
 \end{cases}
\]
For example, for $p\le1$ the potential tends to infinity at both ends
of $0<h<\infty$. Thus all these positive-support solutions stay
between positive finite bounds, continue through all angles, and have
bounded curvature. For $p>1$, if a branch approaches $h=0$ with $E>0$,
then $|h'|\to\sqrt{2E}$ and $|\gamma|\to\sqrt{2E}$.
Its endpoint is a finite arclength away since $ds/dh$ in absolute value
is asymptotic to $h^{p-2}/(C\sqrt{2E})$. Its curvature tends to zero
for $1<p<2$, equals $C$ for $p=2$, and diverges for $p>2$.
Curvature divergence alone does not determine a cusp type; the cusp
calculation above uses the local parametrization.

\paragraph{Turning during one oscillation.}
Suppose $h$ is positive, nonconstant and periodic. Let $P$ denote
its least period in the normal angle. During one period the tangent makes
$\rho=P/(2\pi)$ of a full turn. The position
$\gamma=h e+h'e'$ and the oriented tangent repeat after $n$ oscillations
if $nP=2\pi m$. Conversely, a period of the oriented curve must
advance the normal by an integer number of full turns and must repeat
$h$. Thus the curve closes exactly when $\rho$ is rational. Written in
lowest terms, $\rho=m/n$ means that the first complete traversal has
$n$ oscillations and $m$ tangent turns. When $\rho$ is irrational, each
successive oscillation traces a rotated copy of the preceding arc and
there is no closed traversal.

\begin{corollary}[Which turning fractions occur]
\label{res:9.1}
Fix $C>0$. If $p<3$ and $p\notin\{-1,2\}$, the set of turning fractions
for positive nonconstant periodic support is the open interval whose
endpoints are
\[
 a_p=\frac1{\sqrt{3-p}},\qquad
 b_p=\begin{cases}1/2,& p\le1,\\ 1/(3-p),&1<p<3.\end{cases}
\]
For each reduced rational in this interval there is exactly one such
closed curve up to rotation about $O$ and choice of starting point.
For $p=-1$ the only turning fraction is $1/2$, and for $p=2$ it is $1$.
For $p\ge3$ no positive nonconstant periodic support exists.
\end{corollary}
\begin{proof}
In the force law used by Castelli \cite{castelli} and Rojas \cite{rojas}, the exponent is
$\alpha=p-1$. The two successive extreme values of $r=1/h$ are
separated by angle $P/2$, so $\rho$ equals the apsidal angle divided by
$\pi$. In \cite[Section 3]{castelli} the limits of this angle at zero eccentricity
and at eccentricity tending to one give precisely $a_p$ and $b_p$
after division by $\pi$. The angle is continuous and strictly monotone
between these limits by \cite[Theorem 1.1]{castelli} and \cite[Theorem 1.1]{rojas}, except at
$\alpha=-2,1$. For fixed $C$ the ratio $k=h_{\min}/h_{\max}$ fixes
the two turning points: equating their values of $V_p$ gives
\[
 h_{\max}^{3-p}=\frac{2(1-k^{p-1})}{C(p-1)(1-k^2)}
 \quad(p\ne1),\qquad
 h_{\max}^2=\frac{-2\log k}{C(1-k^2)}\quad(p=1).
\]
Conversely each $0<k<1$ gives two turning points about the unique minimum
of $V_p$, hence a positive periodic solution. The corresponding eccentricity
is $(1-k)/(1+k)$. Uniqueness for the support equation at a turning
point gives uniqueness up to a shift of the angle, hence a rotation.

At $p=-1$ Proposition~\ref{res:4.1} gives centred ellipses; every noncircular
one has least support period $\pi$. At $p=2$, $h''+h=1/C$, and every
nonconstant positive solution has least period $2\pi$; these profiles
are circles with $O$ inside and away from the centre.
For $p\ge3$, the equation at a maximum and a minimum gives
$h_{\max}^{p-3}\le C\le h_{\min}^{p-3}$. For $p>3$ this forces
constancy; for $p=3$ it gives $C=1$ and $h''=0$.
\end{proof}

\section{General relations \texorpdfstring{$B=\Phi(A)$}{B = Phi(A)}}

Locally, a relation between the two densities can be written as
$B=\Phi(A)$. Reversing the orientation replaces $(A,B)$ by $(-A,-B)$, so
the reversed arc satisfies $B=-\Phi(-A)$. A relation of this form therefore
refers to an oriented arc. Except in Example~\ref{res:ex10.1} and
Theorem~\ref{res:10.5}, which concern relations containing $A'$ or the
arclength and have their own hypotheses, we consider in this section
oriented arcs with positive curvature that satisfy $B=\Phi(A)$, where $\Phi$ is a
positive locally Lipschitz function on the range of $A$ under
consideration. Then $\kappa=A^2\Phi(A)>0$,
and in the normal angle the signed support $h=A$ satisfies
\begin{equation}
 h''+h=\frac{1}{h^2\Phi(h)}. \tag{10.1}
\end{equation}
The sign of $AB$, which equals the sign of $\kappa A$, does not depend on
the orientation. For these arcs the case $AB>0$ is positive
support and the case $AB<0$ is negative support.

\begin{remark}[The force law]
\label{res:10.1}
Compare (10.1) with Binet's equation. On an arc with $AB>0$ the polar
reciprocal $r=1/h$ is an orbit of the attractive central force
\[
 F(r)=\frac{L^2}{\Phi(1/r)} .
\]
On an arc with $AB<0$ put $g=-h>0$. Then $g''+g=-1/(g^2\Phi(-g))$,
which is Binet's equation for the repulsive central force of magnitude
$L^2/\Phi(-1/r)$. Thus the angular density, as a function of the areal
density, is $L^2$ divided by the force, and negative support corresponds
to repulsion. The monomial relation $A^pB=C$ is the power law of
Section~9.

For example, let $\mu>0$ and $0\le\nu<1$. The polynomial relation
$A^2B(\mu+\nu A)=1$ is $\kappa(\mu+\nu A)=1$, and (10.1) becomes the linear
equation $h''+(1-\nu)h=\mu$. It corresponds to the sum of an inverse-square
and an inverse-cube force. Its solutions are
$h=\mu/(1-\nu)+a\cos(\sqrt{1-\nu}\,(\theta-\theta_0))$, and they are
positive for all $\theta$ exactly when $|a|<\mu/(1-\nu)$. For
$0<|a|<\mu/(1-\nu)$ the turning fraction is $1/\sqrt{1-\nu}$,
so these curves close exactly when $\sqrt{1-\nu}$ is rational. For
$|a|\ge\mu/(1-\nu)$ the support vanishes at some $\theta$, and the arc with
positive support is not closed.

In particular, if $\sqrt{1-\nu}$ is rational, every solution whose
nonconstant support is positive for all normal angles gives a closed
curve. For $\mu=1$ and $\nu=3/4$ the
solutions are $h=4+a\cos((\theta-\theta_0)/2)$ with $0<|a|<4$, and every
such curve closes after two tangent turns. Thus $A^2B=C$ (circles) and
$B/A=C$ (centred ellipses) are not the only relations all of whose curves
with periodic positive support are closed. Bertrand's theorem
\cite[\S8]{arnold} does not
apply here: it fixes the force and varies the angular momentum, whereas a
fixed relation $B=\Phi(A)$ fixes the quotient $F/L^2$.
\end{remark}

\begin{proposition}[Ends of the arcs with $AB<0$]
\label{res:10.2}
Let $\Phi$ be positive and locally Lipschitz on $(-\infty,0)$, and let
$\gamma$ be a maximal arc with $B=\Phi(A)$ and $A<0$. Put
$\psi(g)=1/(g^2\Phi(-g))$. Then the tangent of $\gamma$ turns through less
than $\pi$, the arc is embedded, and $A\to0$ at each end. Moreover:
\begin{itemize}
\item if $\int_0\psi(g)\,dg$ diverges, each end has infinite length, goes to
infinity, and approaches a ray from $O$;
\item if $\int_0\psi(g)\,dg$ converges, each end has finite length and stops
at a finite point whose tangent line passes through $O$.
\end{itemize}
\end{proposition}
\begin{proof}
Let $g=-h>0$ on the maximal interval $I=(\theta_0,\theta_1)$ of the normal
angle. By (10.1), $g''+g=-\psi(g)<0$, so $g$ is strictly concave. A positive
concave function stays bounded on bounded intervals, so a finite end of $I$
can occur only where $g\to0$.

Consider the end $\theta_1$, assumed finite, and put
$\Delta=\theta_1-\theta$. Near this end $g'<0$, and $|g'|$ increases with
$\theta$ because $g''<0$. Hence
\begin{equation}
 \Delta=\int_0^{g}\frac{dg}{|g'|}\le\frac{g}{|g'(\theta)|},
 \qquad\text{that is,}\qquad |g'|\,\Delta\le g. \tag{10.2}
\end{equation}
The same bound holds at $\theta_0$.

If $I$ contained an interval $[a,a+\pi]$, multiplying $g''+g=-\psi(g)$ by
$w=\sin(\theta-a)$ and integrating would give
$g(a)+g(a+\pi)=-\int\psi(g)\,w\,d\theta<0$, which is impossible. Hence both
ends are finite and $|I|\le\pi$. If $|I|=\pi$, the same computation on
$(\theta_0+\varepsilon,\theta_1-\varepsilon)$ with $a=\theta_0$ has boundary
terms bounded by $g+|g'|\Delta\le2g$ by (10.2). They tend to zero, while the
integral stays negative. Thus $|I|<\pi$. The velocity is
$\gamma'=\psi(g)\,e'(\theta)$, and its projection on $e'$ at the midpoint of
$I$ is positive. So $\gamma$ is a graph over that direction, and it is
embedded.

Fix $g_0$ in the range of $g$. The first integral is
\[
 \tfrac12g'^2+\tfrac12g^2-\int_g^{g_0}\psi(v)\,dv=E,
 \qquad |\gamma|^2=g^2+g'^2 .
\]
Therefore $|g'|\to\infty$ at an end exactly when $\int_0\psi$ diverges.
Let $e_*=e(\theta_1)$ and $e_*'=e'(\theta_1)$. Since $\gamma=-ge-g'e'$,
\[
 \langle\gamma,e_*\rangle=-g\cos\Delta-g'\sin\Delta,\qquad
 \langle\gamma,e_*'\rangle=g\sin\Delta-g'\cos\Delta .
\]
By (10.2) the first is at most $2g$ in absolute value and tends to zero.
The second tends to $\ell=\lim|g'|$. If $\ell=\infty$, the end goes to
infinity along the ray spanned by $e_*'$, and its distance from this ray
tends to zero. If $\ell<\infty$, then $\gamma\to\ell e_*'$. The limiting
tangent direction is $e_*'$, so the tangent line at the endpoint passes
through $O$. The length of the end is
$\int\psi(g)\,d\theta=\int\psi(g)\,dg/|g'|$, which is finite exactly when
$\int_0\psi$ converges.
\end{proof}

For a real exponent $p$ the power $A^p$ need not be real when $A<0$, so on
arcs with $AB<0$ we write the monomial relation as $|A|^pB=C$ with $C>0$;
for an integer $p$ this is $A^pB=(-1)^pC$. Then $\psi(g)=g^{p-2}/C$.
The ends go to infinity for $p\le1$ and stop at a finite point for $p>1$.
The cases $p=0$ and $p=1$ are Proposition~\ref{res:6.3} and the expanding profiles of
Section~5. At a finite endpoint the curvature is the limit of $1/\psi(g)$.

\begin{proposition}[Double points of the closed curves]
\label{res:10.3}
Let $\Phi$ be positive and locally Lipschitz on $(0,\infty)$. Let
$\gamma$ be a closed curve with $B=\Phi(A)$, with positive nonconstant
support, and with turning fraction $m/n$ in lowest terms, traversed once.
Then all self-intersections of $\gamma$ are transversal, no three branches
of $\gamma$ pass through one point, and $\gamma$ has exactly $n(m-1)$
double points.
\end{proposition}
\begin{proof}
Let $\varphi$ be the polar angle of $\gamma$ and $r=|\gamma|$. Since
$r^2\,d\varphi/ds=A>0$, the angle $\varphi$ increases along the curve.
In the normal angle $\theta$ the curve is $\gamma=h\,e(\theta)+h'e'(\theta)$,
so $\varphi=\theta+\arctan(h'/h)$. Since $h>0$, the second term returns to
its initial value after a closed traversal, while $\theta$ advances by
$2\pi m$. Hence $\gamma$ is given by $r=R(\varphi)$ with $\varphi$ running
over an interval of length $2\pi m$. Multiplying (10.1) by $h'$ gives
$r^2=h^2+h'^2=2E-2U(h)$ with $U'(h)=-1/(h^2\Phi(h))<0$, so $r$ increases
with $h$. One oscillation of $h$ rotates the curve by $P=2\pi m/n$, and $h$
is even about each of its extrema. Therefore $R$ is $P$-periodic, even about
its extrema, and strictly monotone between consecutive extrema. Measure
$\varphi$ from a maximum of $R$. Then $R(x)=R(y)$ holds only if
$y\equiv\pm x \pmod P$.

A self-intersection is a solution of $R(\varphi)=R(\varphi+2\pi j)$ with
$1\le j\le m-1$. Since $m$ and $n$ are coprime, the shift
$2\pi j=(jn/m)P$ is not a multiple of $P$. The solutions are therefore
$\varphi\equiv-\pi j \pmod{P/2}$, two in each period and $2n$ in all. Such a
$\varphi$ is not an extremum of $R$, and $R'(\varphi+2\pi j)=-R'(\varphi)\ne0$, so the
two branches cross transversally. Summing over $j$ gives $2n(m-1)$ pairs
$(\varphi,j)$. The pairs $(\varphi,j)$ and $(\varphi+2\pi j,\,m-j)$ describe the same
crossing, which leaves $n(m-1)$. These crossings are distinct points. If
three branches met at one point, three parameter values that differ by
multiples of $2\pi$, but not by multiples of $P$, would give the same value
of $R$. This is impossible, because the values $x$ with a given $R(x)$ form
at most two classes modulo $P$.
\end{proof}

The count $n(m-1)$ was also verified numerically. For
$p=-2,\,0,\,1,\,1.5,\,2.5$ and $21$ turning fractions the closed curves with
$A^pB=1$ were reconstructed and their self-intersections counted directly.

\begin{remark}[Laws with a prescribed turning fraction]
\label{res:10.4}
Write (10.1) as $h''=-V'(h)$ with $V'(h)=h-1/(h^2\Phi(h))$, so that
$E=\tfrac12h'^2+V(h)$. Consider the periodic solutions around a minimum of
$V$, where $V$ is normalized to be zero. For $v>0$ let
$h_-(v)<h_+(v)$ be the two solutions of $V(h)=v$ on either side of the
minimum, and let $D(v)=h_+(v)-h_-(v)$ be the width of the well. Then
\[
 P(E)=\sqrt2\int_0^E\frac{D'(v)}{\sqrt{E-v}}\,dv,\qquad
 D(v)=\frac{1}{\pi\sqrt2}\int_0^v\frac{P(E)}{\sqrt{v-E}}\,dE .
\]
The first formula is the period integral with $v=V(h)$ as the variable on
each side of the minimum. The second is Abel's inversion
\cite[\S12]{landau}; see also \cite{asorey}. Thus the period function
determines the width of the well and conversely, while the midpoint
$M=(h_++h_-)/2$ remains free. Consequently a relation with a prescribed
period function, in particular a strictly monotone one, is constructed as
follows. Compute $D$ from $P$, choose $M$ so that $h_\pm=M\pm D/2$ are
positive and monotone and so that the function
\[
 h(s)=M(s^2)+\tfrac12\operatorname{sgn}(s)\,D(s^2)
\]
extends smoothly through $s=0$ with $h'(s)>0$ on the whole interval of
reconstruction, in particular at $s=0$; this makes the glued
potential, defined by $V(h(s))=s^2$, smooth with a nondegenerate minimum.
Then put
\[
 \Phi(h)=\frac{1}{h^2\bigl(h-V'(h)\bigr)},
\]
which requires $V'(h)<h$. Every smooth relation whose single well has a
nondegenerate minimum arises in this way, subject to these positivity and
invertibility conditions. A degenerate minimum is not covered: for
$V(h)=(h-1)^4$, that is $F(h)=h-4(h-1)^3$, the support function still
oscillates around $h=1$, but $s=\operatorname{sgn}(h-1)\sqrt{V(h)}$ has no
smooth inverse at $h=1$. For example, $P(E)=P_0+\lambda E$ with $P_0>0$ gives
$D(v)=\bigl(2P_0\sqrt v+\tfrac43\lambda v^{3/2}\bigr)/(\pi\sqrt2)$. With a
constant midpoint $M=h_0>0$ this defines, for small energies, a smooth
potential with $V'(h_0)=0<h_0$, and hence a relation whose turning fraction
$P(E)/(2\pi)$ is strictly increasing for $\lambda>0$ and strictly
decreasing for $\lambda<0$. So strictly monotone turning fractions are not
special to the power laws of Section~9. The construction describes these
relations through $P$ and $M$; it does not give a criterion on $\Phi$
itself.
\end{remark}

\paragraph{A relation containing the derivative of $A$.}
Since $A'=\kappa X$ by (2.1), a relation $G(A,B,A')=0$ can be written
as $\kappa=f(A,X)$ wherever it can be solved locally for $\kappa$.
Equations (2.1) then form an autonomous system
in the $(A,X)$ plane, and the curve is recovered from a solution by one
integration. The simplest such relation can be solved explicitly.

\begin{example}[Areal density linear in arclength]
\label{res:ex10.1}
The standing assumptions on $\Phi$ are not used here; the curvature may
take both signs. Consider the relation $A'=c$ with $c$ constant, so that $A=cs+a$.
For $c=0$, (2.1) gives $\kappa X=0$, and the solutions are straight lines
and circles centred at $O$. Let $c\ne0$. Wherever $X\ne0$ the curvature is
$\kappa=c/X$, and the tangent angle $\theta$ is a regular parameter with
$ds/d\theta=X/c$. By (2.1),
\begin{equation}
 A_\theta=X,\qquad X_\theta=\frac Xc-A,\qquad
 A_{\theta\theta}-\frac1c A_\theta+A=0. \tag{10.3}
\end{equation}
Conversely, for a solution $A$ of (10.3) put
$\gamma=A_\theta T(\theta)-A\,N(\theta)$, as in Proposition~\ref{res:8.1}. Then
$\gamma_\theta=(A_\theta/c)\,T$, so the curve is regular where
$A_\theta\ne0$, its areal density is $A$, and $A'=c$ there. Here each
regular arc is oriented so that its unit tangent is $T$; where
$A_\theta/c<0$ this means that the arc is traversed with decreasing
$\theta$. As in the proof
of Proposition~\ref{res:8.2}, a simple zero of $A_\theta$ is an ordinary cusp.
The characteristic equation is $\lambda^2-\lambda/c+1=0$.
\begin{itemize}
\item If $|c|<1/2$, the roots $\lambda_1,\lambda_2$ are real and
$\lambda_1\lambda_2=1$, and
$A=\alpha e^{\lambda_1\theta}+\beta e^{\lambda_2\theta}$. For
$\alpha\beta=0$ the curve is a logarithmic spiral centred at $O$; if it
meets the rays from $O$ at the angle $\psi$, then $\lambda=\cot\psi$ and
$c=\sin\psi\cos\psi$. For $\alpha\beta>0$ the curve is smooth and has no
cusp. One end spirals into $O$ and the other goes to infinity, and the
ratio $X/A$ tends to $\lambda_1$ at one end and to $\lambda_2$ at the
other. Thus the curve passes from the spiral with angle $\psi$ to the
spiral with angle $\pi/2-\psi$. For $\alpha\beta<0$ the curve has one cusp
and one point with $A=0$.
\item If $|c|=1/2$, then $A=(\alpha+\beta\theta)e^{\theta/(2c)}$. For
$\beta=0$ this is the logarithmic spiral with $\psi=\pm\pi/4$, and for
$\beta\ne0$ the curve has one cusp and one point with $A=0$.
\item If $|c|>1/2$, there are no logarithmic spirals. Put
$\omega=\sqrt{1-1/(4c^2)}$. Then
$A=He^{\theta/(2c)}\cos(\omega(\theta-\theta_0))$. The curve is an infinite
chain of arcs of finite length joined at ordinary cusps, one arc for each
interval of length $\pi/\omega$ in $\theta$. Each arc contains exactly one
point with $A=0$; the curve is smooth there, but its tangent line passes
through $O$ and $B$ is not defined. Since
$A(\theta+\pi/\omega)=-e^{\pi/(2c\omega)}A(\theta)$, the next arc is
obtained from the preceding one by the rotation through $\pi+\pi/\omega$
about $O$ and the dilation with factor $e^{\pi/(2c\omega)}$. The whole
curve is invariant under this similarity.
\end{itemize}
\end{example}

\paragraph{A relation containing the arclength.}
The counterpart of the preceding example for the angular density has a
complete answer. Again the standing assumptions on $\Phi$ are not used:
the curve below has an inflection point, and its curvature takes both
signs.

\begin{theorem}[Angular density linear in arclength]
\label{res:10.5}
Let $\gamma$ be a nondegenerate curve with $B=B_0+cs$, where $c\ne0$, and
put $b=|c|$. Then $\gamma$ extends to all $s\in\mathbb R$ as a nondegenerate
curve, and it has exactly one inflection point, at $s_0=-B_0/c$. The
product $AB$ is negative on one side of $s_0$ and positive on the other.
Let $t=|s-s_0|$ be the arclength measured from the inflection point.
\begin{itemize}
\item On the side with $AB<0$,
\[
 |A|\sim\frac{3}{bt^3},\qquad |\kappa|\sim\frac{9}{bt^5},\qquad
 \gamma=(t+d)\,T_\infty+O(t^{-3})\qquad(t\to\infty)
\]
for a unit vector $T_\infty$ and a constant $d$. Thus this end approaches
the ray from $O$ spanned by $T_\infty$.
\item On the side with $AB>0$ there are constants $0<k_1<k_2$ such that
\[
 k_1t^{-1/3}\le|A|\le|\gamma|\le k_2t^{-1/3}
\]
for all large $t$, and $|\kappa|$ lies between positive multiples of
$t^{1/3}$. This end winds infinitely often around $O$ and tends to $O$,
and its length is infinite.
\end{itemize}
\end{theorem}
\begin{proof}
By (2.1) and $\kappa=A^2B$,
\[
 A'=A^2BX,\qquad X'=1-A^3B,\qquad
 (A^2+X^2)'=2X,\qquad (1/A)'=-BX .
\]
The third equation bounds $|\gamma|$ on every finite interval of the
arclength, and then the fourth shows that $A$ cannot reach zero there.
Hence the solution extends to all $s$, and $A$ keeps its sign. Since
$\kappa=cA^2(s-s_0)$, the curvature has exactly one zero, a simple one, and
$AB$ changes sign there.

Fix one side of $s_0$. Let dots denote derivatives with respect to $t$, and
put $u=|A|$, $r=|\gamma|$ and $Y=\langle\gamma,d\gamma/dt\rangle$, so that
$Y=\pm X$, $r^2=u^2+Y^2$ and $(r^2/2)^{\boldsymbol\cdot}=Y$. In these
variables the equations are
\[
 \dot Y=1+btu^3,\quad (1/u)^{\boldsymbol\cdot}=btY\quad(AB<0),\qquad\qquad
 \dot u=btu^2Y,\quad \dot Y=1-btu^3\quad(AB>0).
\]

\emph{The side $AB<0$.} Here $\dot Y\ge1$, so $Y\ge t-K$ for a constant
$K$. Integrating $(1/u)^{\boldsymbol\cdot}=btY$ gives $u=O(t^{-3})$. Then
$\dot Y=1+O(t^{-8})$, so $Y=t+d+O(t^{-7})$ for a constant $d$, and
$1/u=bt^3/3+O(t^2)$. Hence $u\sim3/(bt^3)$ and
$|\kappa|=btu^2\sim9/(bt^5)$. Since $|\kappa|$ is integrable, the unit
tangent $d\gamma/dt$ tends to a unit vector $T_\infty$ with error
$O(t^{-4})$. The component of $\gamma$ along $d\gamma/dt$ is $Y$, and the
component orthogonal to it has absolute value $u=O(t^{-3})$. Therefore
$\gamma=(t+d)T_\infty+O(t^{-3})$.

\emph{The side $AB>0$.} Put
\[
 W=\frac{1}{btu},\qquad H=\tfrac12r^2+W .
\]
Direct differentiation gives
\[
 \dot W=-\frac Wt-Y,\qquad \dot H=-\frac Wt,\qquad
 (YW)^{\boldsymbol\cdot}=3W-2H-\frac{YW}{t}.
\]
For $t\ge t_0>0$ the function $H$ is positive and decreasing, so $r$, $Y$
and $W$ are bounded. Combining the last two identities,
\[
 \frac{d}{dt}\Bigl(t^{2/3}H+\tfrac13t^{-1/3}YW\Bigr)
 =-\tfrac49\,t^{-4/3}\,YW .
\]
The right-hand side is absolutely integrable at infinity, and
$t^{-1/3}YW\to0$. Hence $t^{2/3}H$ has a finite limit. Since
\[
 H\ge\tfrac12u^2+\frac{1}{btu}\ge\tfrac32\,(bt)^{-2/3},
\]
this limit is positive. Therefore $r\le\sqrt{2H}=O(t^{-1/3})$ and
$u=1/(btW)\ge1/(btH)\ge k_1t^{-1/3}$ for some $k_1>0$. Together with
$u\le r$ this gives the two-sided bounds, and $|\kappa|=btu^2$ lies between
positive multiples of $t^{1/3}$. Finally, the polar angle $\varphi$ of
$\gamma$ satisfies $r^2\,d\varphi/ds=A$, so $|d\varphi/dt|=u/r^2$ is bounded
below by a positive multiple of $t^{1/3}$, and $d\varphi/dt$ has constant
sign. Hence the curve winds infinitely often around $O$, while $r\to0$ as
$t\to\infty$.
\end{proof}

The two ends are of the same kinds as those of the Airy curve in
Section~4: one approaches a ray from $O$, and the other winds into $O$
along an arc of infinite length. The next corollary determines the lower
and upper limits of the scaled distance $t^{1/3}|\gamma|$ on the winding
end. They are distinct when $E_*>\tfrac32$.

\begin{corollary}[Limits on the winding end]
\label{res:10.6}
On the side with $AB>0$ in Theorem~\ref{res:10.5} put
$E_*=\lim_{t\to\infty}(bt)^{2/3}H$. Then $E_*\ge\tfrac32$. Let
$q_-\le q_+$ be the positive roots of
\[
 \tfrac12q^2+\frac1q=E_* .
\]
Then
\[
 \liminf_{t\to\infty}t^{1/3}|\gamma|=\liminf_{t\to\infty}t^{1/3}|A|
 =b^{-1/3}q_-,\qquad
 \limsup_{t\to\infty}t^{1/3}|\gamma|=\limsup_{t\to\infty}t^{1/3}|A|
 =b^{-1/3}q_+ .
\]
\end{corollary}
\begin{proof}
The limit exists and is at least $\tfrac32$ by the proof of
Theorem~\ref{res:10.5}. Put
\[
 q=(bt)^{1/3}u,\qquad z=(bt)^{1/3}Y,\qquad
 \tau=\tfrac34\,b^{1/3}t^{4/3},\qquad
 G(q,z)=\tfrac12(q^2+z^2)+\frac1q .
\]
Then $(bt)^{1/3}r=\sqrt{q^2+z^2}$ and $(bt)^{2/3}H=G(q,z)$, and the
equations for the side $AB>0$ become
\[
 \frac{dq}{d\tau}=q^2z+\frac{q}{4\tau},\qquad
 \frac{dz}{d\tau}=1-q^3+\frac{z}{4\tau}.
\]
By Theorem~\ref{res:10.5} the point $(q,z)$ stays in a compact subset of
the half-plane $q>0$ for large $\tau$. Let $\Omega$ be the set of its
limit points as $\tau\to\infty$. The terms $q/(4\tau)$ and $z/(4\tau)$ tend
to zero. Hence, if $(q,z)(\tau_n)\to p$ with $\tau_n\to\infty$, the shifted
solutions $(q,z)(\tau_n+\cdot)$ converge, uniformly on compact intervals,
to the solution through $p$ of the limiting system
\[
 \frac{dq}{d\tau}=q^2z,\qquad \frac{dz}{d\tau}=1-q^3 .
\]
So $\Omega$ contains the whole orbit of the limiting system through each
of its points. The limiting system conserves $G$, and $G(q,z)\to E_*$
along the solution, so $\Omega$ lies in the level set $G=E_*$. The function
$q^2/2+1/q$ is convex with the single minimum $\tfrac32$ at $q=1$, and
$(1,0)$ is the only equilibrium of the limiting system. Therefore, if
$E_*>\tfrac32$, the level set is a single closed orbit, on which $q$ runs
over $[q_-,q_+]$, and $\Omega$ is this orbit. If $E_*=\tfrac32$, the level
set is the point $(1,0)$ and $q_-=q_+=1$. In both cases the lower and
upper limits of $q$ are $q_-$ and $q_+$. On the level set
$q^2+z^2=2E_*-2/q$, which increases with $q$ and equals $q_\pm^2$ at
$q=q_\pm$. Hence the lower and upper limits of $(bt)^{1/3}r$ are also $q_-$
and $q_+$.
\end{proof}

For example, for $c=1$, $A(s_0)=1$ and $X(s_0)=0$ numerical integration
gives $E_*=1.7842$, hence $q_-=0.6308$ and $q_+=1.4929$. The computed
minimum and maximum of $t^{1/3}|A|$ and of $t^{1/3}|\gamma|$ on
$50\le t\le3000$ agree with these values to four decimal places.

\paragraph{Monotonicity of the turning fraction.}
We return to the relations $B=\Phi(A)$ and ask when the turning fraction
$P(E)/(2\pi)$ is a monotone function of the energy. Put
\[
 F(h)=\frac{1}{h^2\Phi(h)},
\]
the radius of curvature as a function of the areal density. Then (10.1) is
$h''+h=F(h)$, and $V'(h)=h-F(h)$ in the notation of
Remark~\ref{res:10.4}. Let $\Phi$ be smooth, and let $h_0>0$ satisfy
$F(h_0)=h_0$ and $F'(h_0)<1$. Then the circle of radius $h_0$ centred at
$O$ satisfies the relation, and $h_0$ is a nondegenerate minimum of $V$
with $V''(h_0)=\omega^2$, where $\omega^2=1-F'(h_0)$. We consider the
oscillations around $h_0$, normalize $V(h_0)=0$, and use the width $D(v)$
of Remark~\ref{res:10.4}. Put $g(v)=\sqrt v\,D'(v)$.

\begin{proposition}[Tests for a monotone turning fraction]
\label{res:10.7}
In this setting the following hold.
\begin{itemize}
\item[(a)] The period is a mean value of $g$:
\[
 P(E)=\sqrt2\int_0^1\frac{g(Ew)}{\sqrt{w(1-w)}}\,dw .
\]
Hence, if $g$ is strictly increasing on $(0,E_1)$, then so is $P$, and
if $g$ is strictly decreasing there, then so is $P$.
\item[(b)] For every energy $E$ of the family,
\[
 P'(E)=\frac1E\int_{h_-(E)}^{h_+(E)}
 \Bigl(1-\frac{2VV''}{V'^2}\Bigr)\frac{dh}{\sqrt{2(E-V)}} ,
\]
where the first factor is continued by its limit $0$ at $h_0$.
\item[(c)] At the circle,
\[
 \lim_{E\to0}P'(E)=\frac{\pi K}{12\,\omega^7},\qquad
 K=5F''(h_0)^2+3\bigl(1-F'(h_0)\bigr)F'''(h_0).
\]
Hence, if $K>0$, the turning fraction is strictly increasing for small
oscillations, and if $K<0$ it is strictly decreasing. If the turning
fraction is a monotone increasing function of the energy, then $K\ge0$,
and if it is monotone decreasing, then $K\le0$.
\item[(d)] The width is a mean value of the period:
\[
 \frac{D(v)}{\sqrt v}=\frac{1}{\pi\sqrt2}\int_0^1\frac{P(vw)}{\sqrt{1-w}}\,dw .
\]
Hence, if $P$ is monotone increasing on $(0,E_1)$, then so is
$D(v)/\sqrt v$, and if $P$ is monotone decreasing, then so is
$D(v)/\sqrt v$. All periods are equal if and only if $D(v)/\sqrt v$ is
constant.
\end{itemize}
\end{proposition}
\begin{proof}
Put $s=\sqrt v$. Near $h_0$ the function $\sigma(h)$ that equals
$\sqrt{V(h)}$ for $h\ge h_0$ and $-\sqrt{V(h)}$ for $h\le h_0$ is smooth
with $\sigma'(h_0)=\omega/\sqrt2$. Let $\psi$ be the inverse of
$h\mapsto\sigma(h)$, shifted so that $\psi(0)=0$. Then
$h_\pm(v)=h_0+\psi(\pm s)$, so $D(v)=\psi(s)-\psi(-s)$ and
\[
 g(v)=\tfrac12\bigl(\psi'(s)+\psi'(-s)\bigr).
\]
This is a smooth even function of $s$, and therefore a smooth function of
$v$ up to $v=0$.

(a) Substitute $v=Ew$ in the first formula of Remark~\ref{res:10.4}. The
weight $1/\sqrt{w(1-w)}$ is positive, which gives the two statements on
monotonicity.

(b) Differentiating (a) gives
\[
 P'(E)=\sqrt2\int_0^1\frac{g'(Ew)\,w}{\sqrt{w(1-w)}}\,dw
 =\frac{\sqrt2}{E}\int_0^E\frac{\sqrt v\,g'(v)}{\sqrt{E-v}}\,dv .
\]
The value $g(v)$ is the sum of the values of $\sqrt V/|V'|$ at the two
points $h_\pm(v)$, and
\[
 \frac{d}{dh}\frac{\sqrt V}{V'}=\frac{V'^2-2VV''}{2\sqrt V\,V'^2},
 \qquad dv=V'\,dh .
\]
Returning to the variable $h$ on both sides gives the formula. Near
$h_0$, with $x=h-h_0$, one has $V'^2-2VV''=O(x^3)$ and $V'^2\sim\omega^4x^2$,
so the first factor tends to $0$.

(c) Write $\psi(s)=c_1s+c_2s^2+c_3s^3+\cdots$. Then
$g(v)=c_1+3c_3v+O(v^2)$, and by the first formula in the proof of (b),
\[
 \lim_{E\to0}P'(E)=\sqrt2\,g'(0)\int_0^1\frac{w\,dw}{\sqrt{w(1-w)}}
 =\frac{3\pi}{\sqrt2}\,c_3 .
\]
Since $V''=1-F'$, $V'''=-F''$ and $V''''=-F'''$, the expansion of $V$ at
$h_0$ is $\tfrac12\omega^2x^2-\tfrac16F''x^3-\tfrac1{24}F'''x^4+\cdots$,
with the derivatives of $F$ taken at $h_0$. Comparing coefficients in
$V(h_0+\psi(s))=s^2$ gives
\[
 c_1=\frac{\sqrt2}{\omega},\qquad
 c_2=\frac{F''c_1^2}{6\,\omega^2},\qquad
 c_3=\frac{c_1^3K}{72\,\omega^4},
\]
and hence the limit $\pi K/(12\omega^7)$. The remaining statements follow
because $P'$ is continuous up to $E=0$.

(d) Substitute $E=vw$ in the second formula of Remark~\ref{res:10.4}. If
$D(v)/\sqrt v$ is constant, then $g$ is constant, and $P$ is constant by
(a).
\end{proof}

Parts (a), (b) and (d) are known for potential systems in other forms; we
have stated them for the relations $B=\Phi(A)$ and proved them from the
two formulas of Remark~\ref{res:10.4}. Part (b) is the case $k=1$ of
\cite[Theorem~A]{manosas}. The same theorem gives a sufficient condition
at pairs of points: the period is monotone if a certain combination of
the values of $1-2VV''/V'^2$ at the two points $h_\pm(v)$ of each level
has no zero. This combination is a positive multiple of $g'(v)$, so the
condition is that of (a). The necessary condition in (d) is
\cite[Theorem~2.6 and Remark~2.7]{manosas}, which extends a result of
Opial \cite{opial} for even potentials; the statement on equal periods is
quoted there from \cite{cima}. Part (c) is the classical first term of
the expansion of the period at a centre, written in terms of $F$.

For the monomial relation $A^pB=C$ with $p<3$ one has $F(h)=h^{p-2}/C$ and
\[
 K=\frac{2\,(p-2)(p+1)(p-3)^2}{h_0^2}.
\]
Thus $K$ vanishes exactly for $p=-1$ and $p=2$, where all periods are
equal; it is negative for $-1<p<2$ and positive for $p<-1$ and for
$2<p<3$. This agrees with Corollary~\ref{res:9.1}, where the turning
fraction decreases from $a_p$ to $b_p$ in the first case and increases in
the other two. For a relation that is not monomial, take
$A^2B(3+A^3)=4$, so that $F(h)=(3+h^3)/4$ and $h_0=1$. Then
$\omega^2=\tfrac14$ and $K=\tfrac{99}{8}>0$, so the turning fraction
increases near the unit circle. 

The condition in (a) is sufficient but not necessary. For $B=1$ one has
$V(h)=\tfrac12h^2+1/h-\tfrac32$, and the period decreases for all energies
by Proposition~\ref{res:7.1}. But the computed values of $g$ decrease only
from $g(0)=\sqrt{2/3}$ to a minimum $0.6955$ near $v=12.2$; after that
$g$ increases towards $1/\sqrt2$, the value that corresponds to the limit
$\pi$ of the period.

\paragraph{No pointwise criterion.}
Proposition~\ref{res:10.7}(b) decides the sign of $P'(E)$ by an integral
over the well, and parts (a) and (d) are conditions on the width $D$,
which involves the two sides of the well at once. Part (c) is a condition
at the single point $h_0$, but it decides only small oscillations. The
next proposition shows that, for smooth relations, no universal condition
on finitely many derivatives of $F$, imposed separately at each point of
the well, characterizes an everywhere positive, or an everywhere negative,
derivative $P'(E)$. It says nothing about weak monotonicity: the comparison
relation in the proof has constant period, which is weakly monotone. The
reason is that a relation with a strictly increasing
turning fraction can coincide with a relation with constant turning
fraction on a whole side of the well.

Consider triples $(F,h_0,E_1)$ of the following kind: $F$ is a smooth
positive function on an open interval $I\subset(0,\infty)$, $h_0\in I$
satisfies $F(h_0)=h_0$ and $F'(h_0)<1$, and the oscillations of the
support function of the curves with
\[
 \frac1\kappa=F(A),\qquad\text{that is,}\qquad B=\Phi(A)\ \text{ with }\
 \Phi(A)=\frac{1}{A^2F(A)},
\]
around $h_0$ with energies in $(0,E_1)$ stay inside $I$. Their union is the
interval $W=(h_-(E_1),h_+(E_1))$, the well. Let $P(E)$ be the period of the
support function, so that $P(E)/(2\pi)$ is the turning fraction.

\begin{proposition}[No pointwise criterion]
\label{res:10.8}
There are no integer $k\ge0$ and function $Q(h,y_0,\dots,y_k)$ such that,
for every triple $(F,h_0,E_1)$ as above, $P'(E)>0$ for all $E\in(0,E_1)$ if
and only if
\[
 Q\bigl(h,F(h),F'(h),\dots,F^{(k)}(h)\bigr)>0
 \quad\text{at every } h\in W\setminus\{h_0\}.
\]
The same holds with $P'(E)<0$ in place of $P'(E)>0$.
\end{proposition}
\begin{proof}
Let $h_0=1$ and $F_0\equiv1$. Then $V_0'(h)=h-1$, so
$V_0(h)=\tfrac12(h-1)^2$, equation (10.1) is $h''+h=1$, and $P_0(E)=2\pi$
for all $E$.

Fix $\varepsilon>0$ and put
\[
 \chi(s)=s\quad(s\ge0),\qquad
 \chi(s)=s\bigl(1+\varepsilon e^{-2/s^2}\bigr)\quad(s<0).
\]
The function $se^{-2/s^2}$, extended by zero, is smooth and flat at $s=0$,
so $\chi$ is smooth with $\chi'(0)=1$. Moreover
$\chi'(s)=1+\varepsilon e^{-2/s^2}(1+4/s^2)>0$ for $s<0$. Hence
$h=1+\chi(s)$ is a smooth increasing change of variable, and
$V_L(1+\chi(s))=\tfrac12s^2$ defines a smooth potential with a
nondegenerate minimum at $h=1$ and $V_L''(1)=1$. Put $F_L(h)=h-V_L'(h)$.
For $h\ge1$ we have $V_L=V_0$ and $F_L=1$. Since $F_L(1)=1$ and
$F_L'(1)=1-V_L''(1)=0$, there is $s_1\in(0,1)$ such that $1+\chi(s)>0$
and $F_L(1+\chi(s))>0$ for $|s|\le s_1$. Take $E_1=s_1^2/2<2$ and let $I$ be
an open interval containing $[1+\chi(-s_1),1+s_1]$ on which $F_L$ is
positive. Then $(F_L,1,E_1)$ is a triple as above, and its well is
$W_L=(1+\chi(-s_1),1+s_1)$.

With $s=\pm\sqrt{2v}$ the two sides of the well are
$h_+(v)=1+\sqrt{2v}$ and $h_-(v)=1-\sqrt{2v}\,(1+\varepsilon e^{-1/v})$, so
\[
 D(v)=\sqrt{2v}\,\bigl(2+\varepsilon e^{-1/v}\bigr),\qquad
 g(v)=\sqrt v\,D'(v)=\sqrt2+\varepsilon\sqrt2\,e^{-1/v}
 \Bigl(\frac12+\frac1v\Bigr),
\]
and
\[
 g'(v)=\frac{\varepsilon\sqrt2\,e^{-1/v}(2-v)}{2v^3}>0
 \qquad(0<v<2).
\]
By Proposition~\ref{res:10.7}(a), $P_L$ is strictly increasing on
$(0,E_1)$, and by the formula in the proof of
Proposition~\ref{res:10.7}(b), $P_L'(E)>0$ for every $E\in(0,E_1)$.

Let $V_R(h)=V_L(2-h)$ and $F_R(h)=h-V_R'(h)=h+V_L'(2-h)$. The substitution
$h\mapsto2-h$ maps the solutions of $h''=-V_L'(h)$ onto those of
$h''=-V_R'(h)$, so $P_R=P_L$. For $h\le1$ we have $F_R(h)=h+(1-h)=1$, and
$F_R(1)=1$, $F_R'(1)=0$. After decreasing $s_1$ if necessary, $F_R$ is
positive near the closure of its well $W_R=(1-s_1,\,1-\chi(-s_1))$, and
$(F_R,1,E_1)$ is a triple as above.

Suppose that $k$ and $Q$ exist. Applied to $F_L$, the condition gives
$Q(h,1,0,\dots,0)>0$ for $1<h<1+s_1$, because $F_L\equiv1$ there. Applied
to $F_R$, it gives the same inequality for $1-s_1<h<1$. Now consider the
triple $(F_0,1,E_1)$ on the interval $(1-s_1,1+s_1)$, whose well is this
whole interval. The condition holds at every point $h\ne1$ of it, so
$P_0'(E)>0$, which contradicts $P_0\equiv2\pi$.

For $P'(E)<0$ take $-\tfrac12<\varepsilon<0$. Then
$\chi'(s)\ge1-\tfrac12\sup_{u>0}e^{-2u}(1+4u)=1-e^{-1/2}>0$, the same
construction applies, and now $g'(v)<0$ on $(0,2)$, so $P_L$ and $P_R$ are
strictly decreasing, and the same argument gives a contradiction.
\end{proof}

The potential $V_L$ was also checked numerically (for $V_R$ the periods are the same): for
$\varepsilon=1$ and $\varepsilon=-\tfrac12$ the periods computed by
integrating $h''=-V_L'(h)$ agree with the formula of
Remark~\ref{res:10.4} to $10^{-10}$, and they increase, respectively
decrease, with $E$.

The relations $F_L$ and $F_R$ are smooth but not real analytic, since each
coincides with the constant $1$ on an interval without being constant. The
proof therefore does not decide whether a pointwise criterion exists for
real-analytic relations, for example for polynomial or rational
relations $B=\Phi(A)$. Proposition~\ref{res:10.8} also leaves room for
conditions at pairs of points: the condition of
Proposition~\ref{res:10.7}(a) and those of \cite[Theorem~A]{manosas}
compare the two points $h_\pm(v)$ of each level, and these are sufficient
for all energies.

\paragraph{What is known.}
\emph{Complete answers for special families.} For $F(h)=h^{p-2}/C$, that
is, for $A^pB=C$, the turning fraction is strictly monotone for every
$p<3$ other than $-1$ and $2$, by the theorems of Castelli
\cite{castelli} and Rojas \cite{rojas} used in Section~9. For
$F(h)=\mu+\nu h$ all periods are equal (Remark~\ref{res:10.1}).
Remark~\ref{res:10.4} describes the general answer from the other side:
within the class of smooth wells with a nondegenerate minimum, an
admissible pair $(P,M)$, consisting of the period function and the
midpoint of the well and subject to the positivity and invertibility
conditions of Remark~\ref{res:10.4}, determines the relation on the well
(not outside it), so the relations with a monotone turning
fraction are those that come from pairs $(P,M)$ with $P$ monotone.

\emph{In the literature.} Since $h''=-V'(h)$ with $V'(h)=h-F(h)$, the
question is that of the monotonicity of the period function of a
one-dimensional potential system, for potentials of the form
$V'(h)=h-F(h)$. For that problem many sufficient conditions are
known; see Opial \cite{opial}, Schaaf \cite{schaaf}, Chicone
\cite{chicone}, and the references in \cite{manosas}. In
\cite[Theorem~A]{manosas} there is a whole sequence of such conditions,
of which Proposition~\ref{res:10.7}(a) is the first. For example, by \cite[Theorem~A]{chicone} the period is
a monotone increasing function of the energy if $V/(V')^2$ is convex on
the well, and a monotone decreasing one if it is concave. Here
monotonicity is meant in the weak sense: for $V=(h-h_0)^2/2$ the quotient
$V/(V')^2$ is constant, and so is the period. Schaaf's condition
\cite[Theorem~1]{schaaf} is pointwise. In our notation it reads as
follows. Put
\[
 K(h)=5F''(h)^2+3\bigl(1-F'(h)\bigr)F'''(h),
\]
so that $K(h_0)$ is the constant $K$ of Proposition~\ref{res:10.7}(c). If
$K(h)>0$ at every point of the well with $F'(h)<1$, and
$(h-F(h))F''(h)>0$ at every point of the well with $F'(h)=1$, then
$P'(E)>0$ for all energies of the family. Thus the quantity that decides
the monotonicity near the circle, together with the additional condition
at the points where $F'=1$, gives a sufficient condition for an increasing
turning fraction for all energies. The additional condition cannot be
dropped. For $F(h)=20-2\sin x$ with $x=h-20$ one has
$V=\tfrac12x^2+2(1-\cos x)$ and $K=20-8\cos^2x+6\cos x\ge6$, and the well
is single. The small oscillations have period $2\pi/\sqrt3$. Since
$V(a)-V(x)=\tfrac12(a^2-x^2)+2(\cos x-\cos a)$, the oscillation of
amplitude $a=2\pi$ is slower than that of $x''+x=0$, and the one of
amplitude $3\pi$ is faster, so the period exceeds $2\pi$ at the first
amplitude and is less than $2\pi$ at the second; it is not monotone. Here
the second condition fails at $x=4\pi/3$. For $A^pB=C$, both hypotheses
of Schaaf's criterion hold throughout the well when $2<p<3$. When $p<-1$
there are no points with $F'=1$, since $F'<0$, but the requirement
$K(h)>0$ fails for sufficiently large $h$, although the turning fraction
is increasing in that case too. By Proposition~\ref{res:10.8}, no pointwise condition of this kind
can also be necessary. A general necessary condition is the one in
Proposition~\ref{res:10.7}(d). When the period is
not monotone, the question becomes how many critical points it has, that
is, energies with $P'(E)=0$. Chow and Sanders \cite{chowsanders} showed
that there are none for cubic potentials $V$ and at most three for
quartic ones; Ma\~nosas and Villadelprat
\cite{manosas} give criteria that bound their number.

\section*{Use of AI tools}

This work was carried out with extensive use of AI assistants, including
Astra by OpenAI and Claude by Anthropic. The author posed the problems,
directed the work and checked the results.

\end{document}